\newcommand{\note}[1]{}
\newcommand{\CR}{\hbox{{$\cal R$}}} 
\newcommand{\CE}{\hbox{{$\cal E$}}}

\newcommand{\CC}{\hbox{{$\cal C$}}}
\newcommand{\CS}{\hbox{{$\cal S$}}}
\newcommand{\CM}{\hbox{{$\cal M$}}}
\newcommand{\CQ}{\hbox{{$\cal Q$}}}

\newcommand{\cg}{\mathfrak{g}}
\newcommand{\cv}{\mathfrak{v}}

\newcommand{\C}{\mathbb{C}}
\newcommand{\Z}{\mathbb{Z}}

\newcommand{\Dsl}{{D\!\!\!\!/}}
\newcommand{\h}{{\scriptstyle\frac{1}{2}}}
\newcommand{\extd}{{\rm d}}
\newcommand{\del}{\partial}

\newcommand{\isom}{{\cong}}
\newcommand{\eps}{{\epsilon}}
\newcommand{\tens}{\mathop{\otimes}}

\newcommand{\Lin}{{\rm Lin}}
\newcommand{\ver}{{\rm ver}}
\newcommand{\Ad}{{\rm Ad}}
\newcommand{\id}{{\rm id}}
\newcommand{\<}{\langle}
\renewcommand{\>}{\rangle}
\newcommand{\End}{{\rm End}}

\renewcommand{\o}{{}_{\scriptscriptstyle(1)}}
\renewcommand{\t}{{}_{\scriptscriptstyle(2)}}
\renewcommand{\th}{{}_{\scriptscriptstyle(3)}}
\newcommand{\fo}{{}_{\scriptscriptstyle(4)}}
\newcommand{\fiv}{{}_{\scriptscriptstyle(5)}}

\newcommand{\bz}{{}^{\bar{\scriptscriptstyle(0)}}}
\newcommand{\bo}{{}^{\bar{\scriptscriptstyle(1)}}}
\newcommand{\bt}{{}^{\bar{\scriptscriptstyle(2)}}}
\newcommand{\uo}{{{}^{\scriptscriptstyle(1)}}}
\newcommand{\ut}{{{}^{\scriptscriptstyle(2)}}}
\newcommand{\umo}{{{}^{\scriptscriptstyle-(1)}}}
\newcommand{\umt}{{{}^{\scriptscriptstyle-(2)}}}

\newcommand{\proof}{{\bf Proof\ }}
\newcommand{\eproof}{$\quad \diamond$}

\newcommand{\und}[1]{{\underline {#1}}}
\newcommand{\eqn}[2]{\begin{equation}#2\label{#1}\end{equation}}

\newcommand{\lbiprod}{{>\!\!\!\triangleleft\kern-.33em\cdot}}
\newcommand{\rbiprod}{{\cdot\kern-.33em\triangleright\!\!\!<}}

\documentclass[11pt]{article}
\usepackage{amssymb,amsmath,epsfig}
\newtheorem{lemma}{Lemma}[section]
\newtheorem{propos}[lemma]{Proposition}

\newtheorem{theorem}[lemma]{Theorem}

\newtheorem{corol}[lemma]{Corollary}
\newtheorem{defin}[lemma]{Definition}

\begin{document}\baselineskip 18pt

{\ }\qquad \hskip 4.3in \vspace{.2in}

\begin{center} {\LARGE RIEMANNIAN GEOMETRY OF QUANTUM GROUPS AND FINITE
GROUPS WITH NONUNIVERSAL DIFFERENTIALS} \\ \baselineskip 13pt{\ }\\
{\ }\\ Shahn Majid \\ {\ }\\ School of Mathematical Sciences,
Queen Mary, University of London\\Mile End Rd, London E1 4NS, UK
\end{center}
\begin{center}
May 2000 -- revised August 2001
\end{center}

\begin{quote}\baselineskip 13pt
\noindent{\bf Abstract} We construct noncommutative `Riemannian
manifold' structures on dual quasitriangular Hopf algebras such as
$\C_q[SU_2]$ with its standard bicovariant differential calculus,
using the quantum frame bundle approach introduced previously. The
metric is provided by the braided-Killing form on the braided-Lie
algebra on the tangent space and the $n$-bein by the Maurer-Cartan
form. We also apply the theory to finite sets and in particular to
finite group function algebras $\C[G]$ with differential calculi
and Killing forms determined by a conjugacy class. The case of the
permutation group $\C[S_3]$ is worked out in full detail and a
unique torsion free and cotorsion free or `Levi-Civita' connection
is obtained with noncommutative Ricci curvature essentially
proportional to the metric (an Einstein space). We also construct
Dirac operators $\Dsl$ in the metric background, including on
finite groups such as $S_3$. In the process we clarify the
construction of connections from gauge fields with nonuniversal
calculi on quantum principal bundles of tensor product form.
\end{quote}

\baselineskip 18pt
\section{Introduction}

Noncommutative geometry has been proposed for many years as a
natural generalisation of geometry to include quantum effects.
Particularly important should be `Riemannian' geometry and moreover
(in our opinion) quantum groups or Hopf algebras should play a
central role\cite{Ma:pla} just as Lie groups do in the classical
case. With such motivation, a systematic formalism of a quantum
groups-based approach to `quantum manifolds' and `quantum
Riemannian manifolds' on  (possibly noncommutative) algebras was
already introduced a few years ago in  \cite{Ma:rie}. We used the
notion of quantum principal bundles (with quantum group fibre) and
connections in \cite{BrzMa:gau}, to define `frame bundle', `spin
connection', `vielbeins' etc. The paper studied both the classical
limit and at the other extreme with the universal differential
calculus (which is formally defined on any algebra). We now follow
up \cite{Ma:rie} with a detailed application of this formalism to
uncover a rich noncommutative Riemannian geometry both of quantum
groups and finite groups equipped with general differential
structures. That q-deformation quantum groups should have a rich
but q-deformed Riemannian geometry is hardly surprising but that we
can encode it, proving as we do in Section~4 that all standard
q-deformations of simple Lie groups are quantum Riemannian
manifolds is a good test of our theory. More surprising perhaps is
that finite groups have as equally rich a  Riemannian geometry as
Lie groups. It is well known that their bicovariant differential
structures are defined by conjugacy classes (this is immediate from
\cite{Wor:dif}), but we now take this much further in Section~5 to
a braided-Lie algebra of invariant vector fields, Levi-Civita spin
connections, Ricci tensor etc. fully analogous to the Lie case. The
formulation of Ricci tensors also make clear that we are in a
position now to do gravitational physics in this noncommutative
setting. In the finite group case functional integration over
moduli spaces of metrics, etc., becomes finite-dimensional
integration. In contrast to lattice approximation the finite
spacing is not an `error' but simply a noncommutative modification
of the geometry which remains exact and hence valid even for a
finite number of points. Meanwhile in the q-deformed case
infinities can be expected to be at least partly regularised as
poles at $q=1$. It may also be\cite{MajSmo:def} that spin-network
quantum gravity in the presence of a cosmological constant should
lead specifically to a q-deformation of conventional Riemannian
geometry. Another application of Hopf algebras to Planck scale
physics is the observable-state duality introduced in\cite{Ma:pla}
and this has been related recently to T-duality in $\sigma$-models
on groups\cite{BegMa:poi}. Also, the first systematic predictions
for astronomical data (for gamma-rays of cosmological origin)
coming out of models with noncommutative spacetime coordinates have
emerged\cite{AmeMa:wav} with measurable effects even if the
noncommutativity is of Planck scale order. In another direction,
noncommutative tori such as studied by Connes, Rieffel and others
have emerged  as relevant to string theory\cite{CDS:non}. Although
we will not attempt such applications here, we do put on the table
a general approach to such models that can be fully computed and
which is (as we show)  adequate to include the rich geometry of
quantum groups and finite groups as basic building blocks, while in
now way limited to them.

{\ }From a mathematical point of view our constructive `bottom up'
approach, in which we build up the layers of geometry more or less
up to (in the present paper) the construction of Dirac operators,
provides a useful complement to the  powerful `top down' approach
of Connes\cite{Con:geo}, in particular, coming out of K-theory and
cyclic cohomology. There one starts with a spectral triple or
`axiomatic Dirac operator' on an algebra as implicitly defining the
noncommutative geometry. It appears that reconciling these two
approaches should be rather important to a full development of both
and this provides a second motivation for the work. Section~5
contains, for example, a first result comparing the constructive
approach with the Connes approach in the case of Dirac operators
built up on finite groups. A physical application of such an
understanding would be in the Connes-Lott approach to the standard
model\cite{Con:rea} where a discrete Dirac operator encodes the
fermion mass matrix. A geometrical way to build up such a $\Dsl$
would translate directly into a prediction for this. A first step
in this direction is in \cite{MaSch:lat}.

An outline of the paper is the following. We recall briefly in
Section~2 the global theory from \cite{Ma:rie}, with general
differential calculi on the base $M$, fibre $H$ and `total space'
algebra $P$ of the (frame) bundle. The new results begin in
Section~3 where we specialise to the `local' theory (the
parallelizable case) where $P=M\tens H$. Most of the work in this
section goes into constructing a suitable nonuniversal differential
structure $\Omega(P)$ and showing that local data such as $V$-bein
and `gauge field' indeed provide a global bundle with soldering
form and global connection. This situation is unusual in that the
global theory is known but until now the trivial bundle theory has
not been constructed as a case of this (other than with the
universal differential structure). What we achieve in this way is a
theory that works at the level of a general algebra $M$ equipped
with a suitable parallelizable differential structure and
associated framing, which is roughly the level of generality that
we are used to in quantum theory by the time one has added
$*$-structures and Hilbert spaces (we do not do this here since we
have enough to do at the algebraic level). It is therefore also the
level of generality appropriate to a definitive `quantum Riemannian
geometry'. Note that a quantum group here is not an essential input
and one could in principle use a more general `coalgebra
bundle'\cite{BrzMa:geo}. The quantum-mechanical meaning of
coalgebra bundles is discussed in \cite{Ma:con}, which also
announces the present results.

In Section~4 we apply this theory to the case where the base $M$
is itself a quantum group. The main result is the construction of
Riemannian metrics for general differential calculi from
$\Ad$-invariant bilinear forms on the underlying braided-Lie
algebra\cite{Ma:lie}, which we apply to standard quantum groups
such as $\C_q[SU_2]$. For completeness also consider the other
extreme of usual enveloping algebras $U(\cg)$ as noncommutative
`flat' spacetimes.

Finally in Section~5 we specialise our theory to finite sets and,
in particular, to finite groups. The main results are in
Section~5.3 where we compute everything for the concrete example
of the permutation group $S_3$ with its order 3 conjugacy class.
We are able to explicitly solve the torsion-free and
metric-compatibility (or `cotorsion-free') equations for the
`braided Killing form' metric and obtain a unique `Levi-Civita'
spin connection. We also compute the Ricci tensor and find that
$S_3$ is essentially an Einstein space, and we compute the natural
Dirac operator. The contribution of the gravitational spin
connection to this is absolutely essential for a charge
conjugation operator or symmetric distribution of eigenvalues
about zero and we consider this a good test of the consistency of
our constructive approach.

Let us note that following \cite{Ma:rie} there have been one or
two other constructive attempts at noncommutative Riemannian
geometry for finite sets and finite groups, see
e.g\cite{DimMul:dis}\cite{Cas:gra}. The first of these (as well as
some earlier works on `Levi-Civita connections' on q-deformed
quantum groups and homogeneous spaces, such as \cite{HecSch:lev})
takes a linear connection $\nabla$ point of view and not a frame
bundle and spin connection one (which is essential for us arrive
at a Dirac operator constructively). Meanwhile \cite{Cas:gra},
while speaking of `vielbeins' and `spin connections' does not
actually provide any form of `metric compatibility' between them
and hence cannot be considered as a theory of gravity at all.
Moreover, there is not any actual noncommutative geometry of the
total space and fibre leading for example to any kind of `Lie
algebra' in which the spin connection should take its values.
These are some of the difficult problems solved in our approach.
Moreover, even if one were interested only in finite groups (say),
it is important that our constructions are not ad-hoc to that case
but `functorial' in the sense of being embedded in a single theory
that works for general algebras and with other limits including
classical and q-deformed ones.

\subsection*{Acknowledgements}

The writing of the manuscript was completed while visiting the CPT
at Luminy, Marseilles during the month of May, 2000; I thank my
hosts there for an excellent stay.

\subsection*{Preliminaries} We use the usual notations for Hopf
algebras as in \cite{Ma:book}, over a general ground field $k$.
Thus $\Delta:H\to H\tens H$ is the coproduct on the algebra $H$,
$\eps:H\to k$ the counit and $S:H\to H$ the antipode, which we
assume to be invertible. The right adjoint coaction of $H$ on
itself is $\Ad_R(h)=h\t\tens (Sh\o)h\th$ in the numerical notation
for the output of repeated coproducts (summation understood). Next,
on any algebra $M$ there is a universal differential calculus with
1-forms $\Omega^1M$ given by the kernel of the product map $M\tens
M\to M$ and $\extd m=1\tens m-m\tens 1$. General or `nonuniversal'
$\Omega^1(M)$ are quotients of this by an $M-M$-bimodule $N_M$.
Also the  universal calculus extends to an entire exterior algebra
with $\extd^2=0$ and a general higher order calculus is a quotient
of that by a differential graded ideal\cite{Con:geo}. Equivalently
one can build up the calculus order by order. Thus $\Omega^1(M)$
has a maximal prolongation by Leibniz and $\extd^2=0$, and
$\Omega^2(M)$ can then be specified as a quotient of the degree 2
part of that, etc.

In the case of a Hopf algebra $H$ one can construct\cite{Wor:dif}
the bicovariant $\Omega^1(H)$ equivalently in terms of crossed
modules $\Omega_0\in
\CM{}_H^H$ where $H$ acts and coacts on $\Omega_0$ from the right
in a compatible manner. Then $\Omega^1(H)=H\tens\Omega_0$ with the
tensor product (co)action from the right and the regular (co)action
of $H$ from the left via its (co)product. The universal calculus in
this case corresponds to $\Omega_0=\ker\eps$ and a general calculus
is a a quotient of this by a right ideal $Q_H$ which is invariant
under the right adjoint coaction. Equally well we can write
$\Omega^1(H)=\bar\Omega_0\tens H$ where $\bar\Omega_0\in {}^H_H\CM$
etc. There is a canonical higher order exterior algebra
characterised by $\extd^2=0$ and the additional relations defined
by quotienting by the kernel of $\id-\Psi$, where
\eqn{worskew}{
\Psi(v\tens_H w)=w\tens_H v,\quad v\in \Omega_0,\quad
w\in\bar\Omega_0.}

A quantum principal bundle\cite{BrzMa:gau} over an algebra $M$ with
universal calculus is $(P,H,\Delta_R)$ where $P$ is an algebra, $H$
a Hopf algebra,  $\Delta_R:P\to P\tens H$ a right coaction and
algebra map, with \eqn{Mfixed}{ M=P^H=\{p\in P|\
\Delta_R(p)=p\tens 1\}\subset P} and $P$ is flat as an
$M$-bimodule, and the sequence \eqn{exact}{ 0\to P(\Omega^1M)P\to
\Omega^1P{\buildrel {\rm ver}
\over\longrightarrow} P\tens \ker\eps} is exact, where ${\rm
ver}(p\tens p')=p\Delta_R(p')$. This is equivalent to a
`Hopf-Galois' extension in the theory of Hopf algebras, e.g.
\cite{Sch:pri}, while arising in this  `differential' form in
\cite{BrzMa:gau}.

For a general calculus $\Omega^1(M)$, a bundle means
$(P,H,\Delta_R)$ as before and also a choice of calculus
$\Omega^1(P)$ and $\Omega^1(H)$ with the former is right-covariant
in the sense $\Delta_RN_P\subset N_P\tens H$ and
\eqn{NM}{ N_M=N_P\cap\Omega^1M\subset\Omega^1P,\quad {\rm
ver}(N_P)=P\tens Q_H.} The first condition here states that we
recover \eqn{OmegMP}{\Omega^1(M)=\{md_Pn|\ m,n\in
M\}\subset\Omega^1(P)} as a restriction, while the second ensures
exactness \eqn{exactnonu}{ 0\to P\Omega^1(M)P\to
\Omega^1(P){\buildrel{\rm ver}_{N_P} \over\longrightarrow} P\tens
\Omega_0} by the induced map ${\rm ver}_{N_P}$. This is equivalent
to the formulation in \cite{BrzMa:gau}, as explained in
\cite{BrzMa:dif}.

\section{Framings and Riemannian geometry with nonuniversal calculi}

Here we briefly recall from \cite{Ma:rie} how the basic definitions
of quantum group gauge theory can be extended to frame bundles,
torsion, metric etc., with new emphasis on the case of general
differential calculus that will concerns us. This is the
noncommutative geometrical picture used in the paper. First of all,
if $V$ is a right $H$-comodule we define
\eqn{assoc}{
\CE=(P\tens V)^H,\quad \CE^*={\rm hom}^H(V,P)} to be `associated'
bundles. They are dual in the sense that composition and
multiplication in $P$ gives a pairing $\CE\tens_M \CE^*\to M$ of
$M$-bimodules (or every element of $\CE^*$ induces a left
$M$-module map $\CE\to M$). This is the same as for the universal
calculus. We further assume natural flatness properties so that
$(P\Omega^1(M))^H=\Omega^1(M)$ etc. We will see these in detail for
tensor bundles.

\begin{defin} A frame resolution of $(M,\Omega^1(M))$ is a quantum bundle
$(P,H,\Delta_R,\Omega^1(P),\Omega^1(H))$ over it as above, a right
$H$-comodule $V$ and an equivariant $\theta:V\to P\Omega^1(M)$
such that the induced left $M$-module map by applying $\theta$ and
multiplying in $P$ is an isomorphism $s_\theta:\CE\isom
\Omega^1(M)$.
\end{defin}

This expresses the cotangent bundle as an associated bundle to a
principal bundle, which is the role of framing. The choice of $H$
is far from unique, however, and need not be any kind of analogue
of $GL_n$. Once framed, vector fields are $\Omega^{-1}(M)\isom
\CE^*$ and similarly for their powers. We call this also a `framing
isomorphism' induced by $\theta$. We then define a quantum metric
as an isomorphism $\CE\isom
\CE^*$, i.e. we require nondegeneracy but do not necessarily
impose any symmetry (which would be unnatural in the
noncommutative theory). When $V$ is finite-dimensional note that
$V^*$ is a left $H$-comodule automatically and we can view $\CE^*$
as given by the same construction as for $\CE$ but with a
left-right reversal and $V$ replaced by $V^*$. We define $\bar
H=H^{\rm op}$ (with the opposite product) and $\bar P=P$ as an
algebra but with the left coaction \eqn{DeltaL}{ \Delta_L p\equiv
p\bz\tens p\bo=S^{-1}p\bt\tens p\bo} in terms of the original
right coaction. Then we have a left-handed bundle and a metric is
equivalent to a {\em coframing} with this bundle and
$V^*,\theta^*:V^*\to \Omega^1(M)P$ giving an isomorphism
$\CE^*\isom \Omega^1(M)$ as right $M$-modules. This is the
`self-dual' generalisation of Riemannian geometry as the existence
of a framing and coframing at the same time. The corresponding
metric is \eqn{metricP}{ g=\sum_a
\theta^*(f^a)\tens_P\theta(e_a)\in \Omega^1(M)\tens_M\Omega^1(M)}
where $\{e_a\}$ is a basis of $V$ and $\{f^a\}$ is a dual basis.
Or to avoid explicitly dualising $V$ we can of course work with
$\theta^*\in \Omega^1(M)P\tens V$ and the metric as the
composition with $\theta$ and $\tens_M$, etc.

Finally, a connection on a quantum principal bundle is an
equivariant complement of $P\Omega^1(M)P\subset \Omega^1(P)$. In
concrete terms this is equivalent to a connection form, which is
an equivariant map \eqn{connP}{ \omega:\Omega_0\to
\Omega^1(P),\quad {\rm ver}_{N_P}\circ\omega=1\tens\id} where we
recall that $\Omega_0$ is a right comodule by the adjoint coaction
(as part of the crossed module structure). The associated
projection $\Pi_\omega=\cdot_P(\id\tens\omega){\rm ver}_{N_P}$
defines a covariant derivative \eqn{covP}{ D_\omega:\CE\to
\Omega^1(M)\tens_M \CE,\quad D_\omega=(\id-\Pi_\omega)\circ
\extd\tens \id} provided  $(\id-\Pi_\omega)\circ \extd P\subset
\Omega^1(M)P$, in which case one says that $\omega$ is strong. It
is clear that  a (strong) connection $\omega_U$ on the bundle with
universal calculus such that $\omega_U(Q_H)\subset N_P$ induces
one on the bundle with general calculus. In the presence of a
framing, we define:

\begin{defin} Associated to strong $\omega$ is the covariant derivative
$\nabla_\omega:\Omega^1(M)\to\Omega^1(M)\tens_M\Omega^1(M)$
according to the framing isomorphism $s_\theta$, namely
$\nabla_\omega =(\id\tens s_\theta)\circ D_\omega\circ
s_\theta^{-1}$.
\end{defin}

Both $D_\omega$ and hence $\nabla_\omega$ behave in the expected
way with respect to left-multiplication by $M$. One can then
proceed to identify other geometrical objects in terms of
$\omega,\theta$. Thus, torsion  \eqn{torP}{ T:\Omega^1(M)\to
\Omega^2(M)} corresponds under framing isomorphisms to $\bar
D_\omega\wedge \theta:V\to P\Omega^2(M)$ (here we need a
left-handed version of the bundle as explained in \cite{Ma:rie}.)
Specifically, we apply this in the same manner as the construction
of $s_\theta$ to give a map $\CE\to \Omega^2(M)$ which becomes $T$
as stated under $s_\theta$. In this self-dual formulation it is
natural to ask also that the `cotorsion' vanishes. This is the
torsion of $\omega$ with respect to the coframing, i.e.
$D\theta^*\in \Omega^1(M)\tens_M\CE$ which we view via  $s_\theta$
as
\eqn{cotorP}{ \Gamma\in \Omega^2(M)\tens_M\Omega^1(M).} Its
vanishing is a generalisation of `metric compatability' as
explained in \cite{Ma:rie}. Note that the vanishing torsion and
cotorsion require us to specify $\Omega^2(M)$ suitably. We look at
this in detail for trivial bundles in the next section. Similarly,
the Riemann curvature is
\eqn{Riemann}{R:\Omega^1(M)\to\Omega^2(M)\tens_M\Omega^1(M)} as
left $M$-module map corresponding to the curvature of $\omega$.
With some mild additional structure we can also define the Ricci
tensor by a contraction. The most explicit, which we will adopt, is
to apply lift $i:\Omega^2(M)\to \Omega^1(M)\tens_M\Omega^1(M)$ and
take a trace as an $M$-module map with values in the remaining
$\Omega^1(M)\tens_M\Omega^1(M)$. One could also view this as
associated to an interior product or a Hodge $*$-operation.   Let
us also note that once $\Omega^2(M)$ is specified one could impose
a `symmetry' condition on the metric if desired, as in the kernel
of the wedge product
\eqn{metricsym}{ \wedge(g)=0.}

Finally, we discuss some general aspects in this context of `Dirac
operator'.  Most of the definition is straightforward; we define a
spinor as $\psi\in \CS=(P\tens W)^H$ the associated bundle to some
other representation of $H$. Since $H$ is not required to be
anything like $SO_n$ but can be a more general framing it is not
necessary to speak here of double covers or lifting; we simply
frame by the more suitable quantum group to begin with. Then
$D_\omega\psi\in \Omega^1(M)\tens_M \CS$ maps over under the
framing to $\CE\tens_M\CS$. The missing data to define an operator
$\Dsl:\CS\to \CS$ with reasonable properties under scalar
multiplication of spinors is therefore a left $M$-module map
\eqn{biggamma}{ \und\gamma:\CE\tens_M\CS\to \CS.} Classically,
this would be induced by a map $\gamma:V\tens W\to W$ with
equivariance and `Clifford algebra' properties with respect to the
metric. Note also that in place of an `inner product' on $\CS$ it
is natural in our self-dual formulation to have instead an adjoint
spinor space $\CS^*=\hom^H(W,P)$ and $\Dsl$ defined on this
similarly with $\und\gamma^*$. We do not attempt here a full
formulation but will look at some of these issues for trivial
bundles and quantum groups.

\section{Parallelizable Riemannian structures on algebras}

In this section we apply the formalism above to obtain a general
class of quantum Riemannian manifold structures on algebras $M$ for
which the quantum frame bundle has the tensor product form
$P=M\tens H$, i.e. the parallelizable case. Other trivialisations
can change this form, i.e. we work in what we call the {\em tensor
product gauge}. Our main result is the construction of
$\Omega^1(P)$ such that the global theory above is induced from a
`local' theory where global connections correspond to gauge fields
$A:\Omega_0\to\Omega^1(M)$ and soldering forms to $V$-beins $e:V\to
\Omega^1(M)$. The choice of  $\Omega^1(P)$ is far from obvious, for
example $N_P$ generated as a $P$-bimodule by $N_M,N_H$ as suggested
in \cite{BrzMa:dif} would not allow these correspondences to
proceed.

\begin{propos} On $P=M\tens H$ with $\Omega^1(M),\Omega^1(H)$ given, we take
$\Delta_R$, $\Omega^1(P)$ defined by
\[ \Delta_R=\id\tens\Delta, \quad N_P=N_M\tens H\tens H+M\tens M\tens N_H
+\Omega^1 M\tens \Omega^1H\] where we identify $P\tens P=M\tens
M\tens H\tens H$. Then $(P,\Omega^1(P),\Delta_R)$ is a quantum
principal bundle with nonuniversal calculus over $M,\Omega^1(M)$.
Moreover, we may identify the $H$-comodules \[
\Omega^1(M)P=P\Omega^1(M)P=P\Omega^1(M)=\Omega^1(M)\tens H.\]

\end{propos}
\proof The coaction $\Delta_R$ is only on the $H\tens H$ part and
each component of $N_P$ is clearly invariant under this. Hence
$\Delta_R(N_P)\subset N_P\tens H$. Also \[ \ver=\cdot_M\tens
\ver_H,\quad {\rm ver}(m_i\tens n_i\tens h_i\tens g_i)=m_in_i\tens
h_ig_i\o\tens g_i\t\] for $m_in_i\tens h_ig_i=0$  has
$\ver(Omega^1M\tens H\tens H)=0$ and hence ${\rm ver}(N_P)=P\tens
Q_H$ as required (here ${\rm ver}_H$ corresponds to $H$ as a bundle
over $k$).  Next we note that for any algebras $M,H$,\[ M\tens
M=M\tens 1\oplus \Omega^1 M=1\tens M\oplus
\Omega^1 M,\quad H\tens H=H\tens 1\oplus \Omega^1 H=1\tens H\oplus
\Omega^1 H\] by identifying $m\tens n=mn\tens 1-m\extd n$ or
$m\tens n=1\tens mn-(\extd m)n$ for the two cases and similarly
for $H\tens H$. Hence (making choices, i.e. not canonically) we
can write \[ N_P=N_M\tens 1\tens H\oplus M\tens 1\tens N_H\oplus
\Omega^1M\tens\Omega^1H\] as a vector space. From this it is clear
that $N_P\cap \Omega^1M\tens 1\tens 1=N_M\tens 1\tens 1$ as
required. Hence we have a quantum principal bundle. Also from a
similar decomposition we identify \[ N_P\cap P\Omega^1M=N_M\tens
H\tens 1,\quad N_P\cap (\Omega^1 M)P=N_M\tens 1\tens H\]  and hence
we  can identify  $\Omega^1(M)P=\Omega^1(M)\tens 1\tens H$ and
$P\Omega^1(M)=\Omega^1(M)\tens H\tens 1$. Finally,
\[ N_P\cap P(\Omega^1M)P=N_M\tens 1\tens H\oplus
\Omega^1M\tens\Omega^1H=N_M\tens H\tens 1\oplus
\Omega^1M\tens\Omega^1H\] so that we can identify $P\Omega^1(M)P$
with either $\Omega^1(M)P$ or $P\Omega^1(M)$. When the context is
clear we therefore omit the $\tens 1$ and identify all three with
$\Omega^1(M)\tens H$. It remains to verify that these
identifications are $\Delta_R$-covariant, in particular that of
$P\Omega^1(M)P$. We need for this that the identifications $H\tens
H\isom 1\tens H\oplus \Omega^1H$ etc., are equivariant under the
tensor product of the coaction $\Delta$ in each factor up to an
error in $\Omega^1H$. In particular the projection to $1\tens H$
by multiplication is covariant just because $\Delta$ is an algebra
homomorphism. \eproof

As a justification for this calculus note that classically the
three spaces $P\Omega^1(M)$, $\Omega^1(M)P$ and $P\Omega^1(M)P$
coincide, which we have arranged also here. It means that all
connections are automatically strong, etc, as in the classical
theory. Also, $\Omega^1(P)$ has the right size. Thus, for any (say)
finite-dimensional algebra $M$ define
\eqn{dimM}{ \dim(\Omega^1(M))=\dim(M)-1-\frac{\dim(N_M)}{\dim(M)}.}
which is the dimension over $M$ in the free case. Then for the
above $\Omega^1(P)$ we have
\eqn{dimP}{ \dim(\Omega^1(P))=\dim(\Omega^1(M))+\dim(\Omega^1(H)).}

Next we consider framings and coframings with the above
$\Omega^1(P)$ understood. As for the universal calculus in
\cite{Ma:rie} we define to this end a `$V$-bein' and `$V$-cobein'
as linear maps \eqn{eV}{ e:V\to
\Omega^1(M),\quad e^*:V^*\to
\Omega^1(M)}
such that there are induced isomorphisms \[ s_e:M\tens V\isom
\Omega^1(M),\quad s_{e^*}:V^*\tens M\isom\Omega^1(M),\quad
s_e(m\tens v)=me(v),\quad  s_{e^*}(w\tens m)=e^*(w)m.\]

\begin{propos} A framing and coframing of $M$ with $P=M\tens H$ are
equivalent to $(V,e,e^*)$ where $V$ is a right $H$-comodule and
$e$, $e^*$ are a $V$-bein and $V$-cobein in the sense above. The
(co)frame resolutions and quantum metric are  \[
\theta(v)=e(v\bo)\tens v\bt\tens 1,\quad
\theta^*(w)=e^*(w\bo)\tens 1\tens w\bt,\quad g=\sum_a
e^*(f^a)\tens_M e(e_a).\]
\end{propos} \proof Note first that $(H\tens V)^H\isom V$ by $\eps$
in one direction and conversely by $v\mapsto S^{-1}v\bt\tens v\bo$,
hence $\CE\isom M\tens V$. Likewise $\hom^H(V,H)\isom V^*$ by
composing with $\eps$ in one direction and $w\mapsto \phi(w)$,
$\phi(w)(v)=\<w,v\bo\>v\bt$, hence $\CE^*\isom V^*\tens M$. This
part of the standard analysis for associated bundles in the trivial
case\cite{BrzMa:gau}.  Given $e,e^*$ we define respectively
$\theta,\theta^*$ as stated and verify they are equivariant. Thus
$\theta(v\bo)\tens v\bt=e(v\bo\bo)\tens v\bo\bt\tens v\bt\tens
1=\Delta_R\theta(v)$ as $V$ is a right comodule. Similarly for
$\theta^*$ where $V^*$ is a right $H$-comodule by
$\<v,w\bo\>w\bt=\<v\bo,w\>S^{-1}v\bt$ as usual (i.e. the adjoint of
the left $H$-comodule structure on $V$ corresponding in the manner
of (\ref{DeltaL}) the right comodule structure on $V$). Finally the
induced  \[s_\theta(m\tens h\tens v)=(m\tens h)e(v\bo)\tens
v\bt\tens 1=me(v\bo)\tens hv\bt\tens 1\] under the above
identification becomes  \[m\tens v\mapsto s_\theta(m\tens
S^{-1}v\bt\tens v\bo)=me(v\bo\bo)\tens (S^{-1}v\bt)v\bo\bt\tens
1=me(v)\tens 1\tens 1\] i.e. reduces to $s_e$. Likewise
$s_{\theta^*}$ reduces to $s_{e^*}$. Hence we obtain framings and
coframings respectively from $e,e^*$. Conversely any equivariant
$\theta,\theta^*$ must have this form by similar arguments as for
$\CE,\CE^*$. Given these, the general formula for the metric then
reduces to the one shown on using invariance of $f^a\tens e_a$. In
fact the computation here is the same as for the universal calculus
and works for any reasonable calculus on $P$ where
$\Omega^1(M)P\isom \Omega^1(M)\tens 1\tens H$ etc. For our
particular $\Omega^1(P)$ we can suppress the $\tens 1$ in the
formulae for $\theta,\theta^*$. \eproof

Next, for the principal bundle $P=M\tens H$ a trivial reference
connection is provided by
\eqn{refomega}{ \omega_0(v)=1\tens 1\tens \pi_{N_H}(S\tilde v\o\tens \tilde v\t)}
where $\tilde v\in\ker\eps$ is any lift of $v\in\Omega_0$ and
$\pi_{N_H}$ the projection to $\Omega^1(H)$ (the Maurer-Cartan
form of $H$ viewed in $\Omega^1(P)$). Here we view
$\Omega^1(H)\subset \Omega^1(P)$ by the same arguments as for
$\Omega^1(M)$ (their situation is symmetric). Any other connection
then corresponds to the addition of an $\Ad$-equivariant form in
the kernel of $\ver_{N_P}$, i.e.  $\omega-\omega_0:\Omega_0\to
P\Omega^1(M)P$. For our choice of $\Omega^1(P)$ the target here
can be identified with $\Omega^1(M)P$.

\begin{theorem}  A connection on $\Omega^1(P)$ is equivalent to a
linear map or `gauge field'
\[ A:\Omega_0\to \Omega^1(M).\]
The resulting connection and corresponding projection are
\[ \omega(v)= \omega_0(v)+ \pi_{N_P}((S \tilde v\o)\cdot
A(\pi_{\Omega_0}\tilde v\t)\cdot\tilde v\th)\]
\[ (\id-\Pi_\omega)(m_i\tens n_i\tens h_i\tens g_i)=-m_in_i
A(\pi_{\Omega_0}g_i\o)\tens 1\tens h_i g_i\t +m_i\extd n_i\tens
1\tens h_ig_i\in \Omega^1(M)P\] in a manifestly strong form. Here
$\tilde v$, $m_i\tens n_i\tens h_i\tens g_i$ are representatives in
$\ker\eps$ and $\Omega^1P$ respectively and $\pi_{\Omega_0}$
denotes the canonical projection to $\Omega_0$, etc. \end{theorem}
\proof  For any $H$-comodule $V$ we identify equivariant maps
$V\to \Omega^1(M)P$ with linear maps $V\to \Omega^1(M)P$ by the
same construction as above for $V\to P$. Thus $A:V\to \Omega^1(M)$
corresponds to $\tilde A(v)=A(v\bo)\tens 1\tens v\bt$ and
conversely every $\omega$ has this form. In particular we take
$V=\Omega_0$ and the right adjoint coaction given by projecting
down that on $\ker\eps$. Thus \[ \tilde A(v)=A(\pi_{\Omega_0}\tilde
v\t)\tens 1\tens (S\tilde v\o)\tilde v\th.\] When we identify
$\Omega^1(M)P$ with $P\Omega^1(M)P$ we obtain the form for
$\omega-\omega_0$ shown. Note that
\[ \pi_{N_P}(A(\pi_{\Omega_0}\tilde v\t)\tens S\tilde v\o\tens
\tilde v\th)=
\pi_{N_P}(A(\pi_{\Omega_0}\tilde v\t)\tens 1\tens (S\tilde
v\o)\tilde v\th)\] so that the left hand side is manifestly
well-defined. Here the difference between the expressions is in
$\Omega^1(M)\tens\Omega^1H$ and hence killed by the form of $N_P$.
Conversely it clear that $\omega-\omega_0$ is necessarily of this
form as explained. Finally. given such a connection, we have from
the form of $\ver_{N_P}$, the corresponding projector
\[ \Pi_\omega(m_i\tens n_i\tens h_i\tens g_i)=\pi_{N_P}(m_in_i
A(\pi_{\Omega_0}g_i\th)\tens h_i g_i\o\tens (Sg_i\t)g_i\fo))
+m_in_i\tens 1\tens \pi_{N_H}(h_i\tens g_i)\] for any
representative $m_i\tens n_i\tens h_i\tens g_i\in \Omega^1P$. Under
$\pi_{N_P}$ we can move the $h_ig_i\o$ to the second factor and
cancel using the antipode axioms. We also write $m_in_i\tens
1=m_i\tens n_i-m_i\extd n_i$ and $m_i\extd n_i\tens h_i\tens
g_i=m_i\extd n_i\tens 1\tens h_ig_i$  under $\pi_{N_P}$. In this
form we have no further quotient and drop the $\pi_{N_P}$ as shown.
Note that if $h_i\tens g_i\in N_H$ then $h_ig_i\o\tens g_i\t\in
H\tens Q_H$, but since $Q_H$ is $\Ad$-invariant we have $h_i
g_i\o\tens g_i\th\tens (Sg_i\t)g_i\fo\in H\tens Q_H\tens H$.
Multiplying the two copies of $H$ we conclude that $h_ig_i\t\tens
g_i\o\in H\tens Q_H$ also. Therefore $\id-\Pi_\omega$ is
well-defined. \eproof

Note that we do not consider here the question of gauge
transformations themselves, which is much more subtle for
nonuniversal calculi even when the bundle is trivial: we simply
show that all connections in our `tensor product gauge' have the
above form. Basically, a gauge transformation changes the
description of the bundle to a cocycle cross product as explained
in \cite{Ma:diag}, which in turn changes the description of the
calculus (this is a quantum effect in that one does not have this
cocycle classically). Other trivialisations and correspondingly the
formulae in other gauges can in principle be computed via a bundle
automorphism if one wants formulae for `gauge theory' but the
tensor product form of the bundle $P$ will also transform.

\begin{propos} Given a gauge field on $M$ as above and $V$ any right
$H$-comodule, the vector spaces $E=M\tens V$ and $E^*={\rm
Lin}(V,M)$ acquire covariant derivatives
\begin{eqnarray*}&& D_A:E^*\to \Omega^1(M)\tens_M E^*,
\quad (D_A\sigma)(v)=\extd\sigma(v) -\sigma(v\bo)\cdot
A(\tilde\pi_{\Omega_0} v\bt),\\
&& D_A:E\to \Omega^1(M)\tens_M E,\quad D_A\psi=(\extd\tens\id)\psi
-\psi_iA(\tilde\pi_{\Omega_0}\psi^i\bz)\tens
\psi^i\bo,\end{eqnarray*} where $\psi=\psi_i\tens\psi^i\in M\tens
V$ is a notation and $\tilde\pi_{\Omega_0}$ denotes projection to
$\ker\eps$ followed by $\pi_{\Omega_0}$.
\end{propos}
\proof Given $\sigma\in E^*$ we view it as $\Sigma\in \hom^H(V,P)$
as usual by $\Sigma(v)=\sigma(v\bo)\tens v\bt$. Then
\begin{eqnarray*} (\id-\Pi_\omega)(\extd \Sigma(v))& &
=(\id-\Pi_\omega)[1\tens
\sigma(v\bo)\tens 1\tens v\bt-\sigma(v\bo)\tens 1\tens v\bt\tens 1] \\
& &=\extd\sigma(v\bo)\tens 1\tens
v\bt-\sigma(v\bo)A(\tilde\pi_{\Omega_0} v\bt\o)\tens 1\tens
v\bt\t.\end{eqnarray*} However, this equivariant map $V\to
\Omega^1(M)P$ is in the image of the identification (as in the
proposition above) with $\Lin(V,\Omega^1(M))=\Omega^1(M)\tens_M
E^*$  of $D_A\sigma$ as stated. Similarly for $D_A\psi$. One may
verify directly that both maps are well-defined. \eproof

These formulae are characterised not by gauge covariance but by
the global constructions of the previous section specialised to
the case of a tensor product bundle. They are the basic local
formulae of quantum group gauge theory with nonuniversal calculus
in the tensor product gauge.   Now we suppose the existence of
$V$-(co)beins or framings and coframings as explained above. Then
$D_A$ induces $\nabla_A$ etc. under the framing isomorphisms:

\begin{corol} The covariant derivative $\nabla_A:\Omega^1(M)\to
\Omega^1(M)\tens_M\Omega^1(M)$ is given by
\[ \nabla_A=\extd s^{-1}_{ei}\tens_M e(s_e{}^{-1i})-s^{-1}_{ei}\cdot A(
\tilde\pi_{\Omega_0} s_e{}^{-1i}\bz)\tens_M e(s_e{}^{-1i}\bo),\]
where $s^{-1}_{ei}\tens s_e{}^{-1i}$ denotes the output of
$s_e^{-1}$ and we use the projected right adjoint coaction viewed
as a left coaction as in (\ref{DeltaL}). If we write
$\alpha=\alpha^a\cdot e(e_a)$ for all $\alpha\in\Omega^1(M)$, then
this is \[ \nabla_A\alpha=\extd \alpha^a\tens_M e(e_a)-\alpha^a
A(\tilde\pi_{\Omega_0} e_a\bz)\tens_M e(e_a\bo).\] \end{corol}

Similarly for trivial bundles we can look at the construction of
$\und\gamma$. Here $\CS$ can be identified with $S=M\tens W$ as a
left $M$-module as explained above for any associated bundle.

\begin{corol} For $P=M\tens H$ and given $s_e$ and a right-comodule $W$,
suitable $\und\gamma$ in (\ref{biggamma}) are provided by linear
maps $\gamma:V\to \End(W)$. The corresponding Dirac operator $S\to
S$ is given on $\psi=\psi_i\tens\psi^i\in M\tens W$ by \[
\Dsl\psi= \del^a\psi_i \tens \gamma_a(\psi^i)-\psi_i
A^a(\tilde\pi_{\Omega_0}\psi^i\bz)\tens \gamma_a(\psi^i\bo),\quad
s_e^{-1}\circ \extd=\del^a\tens e_a,\quad \gamma_a=\gamma(e_a)\]
where $A^a$ are the components of $A$ as above.
\end{corol} \proof Since $\und\gamma$ is a left $M$-module map
and defined on $(M\tens V)\tens_M (M\tens W)\isom M\tens V\tens
W$, it is determined by $\und\gamma((1\tens v)\tens_M(1\tens
w))\equiv \und\gamma(1\tens v\tens w)\equiv \gamma(v)(w)\in M\tens
W$, say. It is natural to assume here that $\gamma(v)(w)\in W$
itself. Note that the right $M$-module structure on $M\tens V$ is
not the obvious one (it is the one corresponding to that of
$\Omega^1(M)$ via $s_e$) but becomes irrelevant after we absorb
$\tens_M M$. We then compute $\Dsl$  by the above formulae for
$D_A$ on $S$ and the left $M$-module isomorphism $s_e^{-1}$ as
before (with the notations stated) to map $\extd \psi$ and $A$
over to $M\tens V$, thereby obtaining an element of $M\tens V\tens
W$. We then apply $\gamma$ to $V$ and evaluate its output in
$\End(W)$ on the other (spinor) component of $\psi$. \eproof

We note that the operators $\del^a$ in these expressions are not
derivations but characterised by \eqn{partialM}{
\del^a(mn)=m(\del^a n)+(\del^bm)\rho_b{}^a(n);\quad } where we
write the `generalised braiding' or entwining operator induced by
$s_e$ as \eqn{partialbraid}{ \Psi_e:V\tens M\to M\tens
V,\quad\Psi_e(e_a\tens m)=s_e^{-1}(e(e_a)m)=\rho_a{}^b(m)\tens
e_b} for operators $\rho_a{}^b$ on $M$.  They evidently obey
$\rho_a{}^b(1)=\delta_a{}^b$ and
$\rho_a{}^b(mn)=\rho_a{}^c(m)\rho_c{}^b(n)$ as an expression of
the right module structure of $\Omega^1(M)$. In this notation,
\eqn{DslM}{ [\Dsl,m]=(\del^am)\rho_a{}^b\tens \gamma(e_b)} if one
wants to compare this approach with that of Connes\cite{Con:geo}.
From this it is clear that if $\gamma:V\to \End(W)$ is injective
then $\ker \pi_{\Dsl}=N_M$, where $\pi_{\Dsl}(m\extd
n)=m[\Dsl,n]$. Hence these approaches correspond to the same
differential calculus at degree 1. At higher degree Connes
proposes to quotient the universal exterior algebra by the
differential ideal generated from repeated commutators with
$\Dsl$. At degree 2 the requirement that we recover a given choice
of $\Omega^2(M)$ is a quadratic constraint on the linear maps
$\gamma$ appearing in Corollary~3.6. Another aspect to the
`correct' choice of $\gamma$ would be to demand that it is
$H$-equivariant as an analogue of the idea that the gamma-matrices
generate a representation of the spin group. We will look at these
constraints in detail in the settings of Sections~4 and~5.

We require similar properties as in Proposition~3.1 for
$\Omega^2(M)$ and $\Omega^2(P)$ needed for the global picture of
curvature, torsion and cotorsion. Namely, we require
\eqn{deg2P}{\Omega^2(M)\subset\Omega^2(P),\quad \Omega^2(M)P
\subset\Omega^2(P)}
etc. in the obvious way by $\tens 1$ (as above for 1-forms). For
example $\Omega^1(P)$ itself determines a `maximal prolongation' to
higher forms consisting of $\Omega^1(P)\tens_P\Omega^1(P)$ modulo
the additional relations implied by extending $\extd:\Omega^1(P)\to
\Omega^2(P)$ with a graded Leibniz rule and $\extd^2=0$, and a short
computation shows that this works. A general choice will be a
bimodule quotient of this. Similarly for higher degree. We may then
proceed to make calculations along exactly the same lines as for
1-forms above. Specifically, it is clear that $\Lin(V,\Omega^2(M))$
corresponds in the same manner as before to equivariant maps $V\to
\Omega^2(M)P$, etc. One has therefore \[
D_A:\Lin(V,\Omega^n(M))\to \Lin(V,\Omega^{n+1}(M)),\quad
D_A\sigma(v) =\extd\sigma(v)+(-1)^{n+1}\sigma(v\bo)\wedge
A(\tilde\pi_{\Omega_0}v\bt)\] etc.  Here $D_A$ is $\extd$ on $P$
followed by $(\id-\Pi_\omega)$ in each copy of $\Omega^1(P)$. The
proof is just as for the universal calculus in \cite{BrzMa:gau}
followed by the required projections. See also \cite{Ma:diag}.

\begin{propos} For all $\sigma\in\Lin(V,M)$,
$D_AD_A\sigma(v)=-\sigma(v\bo)F_A(\pi_\eps v\bt)$, where
\[ F_A:\ker\eps\to \Omega^2(M),
\quad F_A(v)=\extd
A(\tilde\pi_{\Omega_0}v)+A(\tilde\pi_{\Omega_0}v\o)\wedge
A(\tilde\pi_{\Omega_0}v\t)\] and $\pi_\eps(h)=h-\eps(h)$. We say
that $A$ is `regular' if $F$ descends to $\Omega_0\to
\Omega^2(M)$, i.e. if
\[ A(\tilde\pi_{\Omega_0}q\o)\wedge A(\tilde\pi_{\Omega_0}q\t)=0,
\quad \forall q\in Q_H.\]
\end{propos}
\proof We apply the above formulae for $D_A$ and compute exactly
as for the universal calculus. As in the usual computation
iteration of the coaction produces a coproduct and the
well-defined formula for $D_AD_A\sigma(v)$ as stated. We omit
details since they as the same as the universal case in
\cite{BrzMa:gau}. See also \cite{Ma:diag}. The map
$A\circ\pi_{\Omega_0}$ plays the role of $A:\ker\eps\to\Omega^1M$ in the
universal calculation and all expressions are finally projected to
the relevant differentials. In doing this one only knows that
$F_A:\ker\eps\to \Omega^2(M)$ as stated. \eproof

It is not such a problem if $A$ is not regular. Classically it
would mean that $F_A$ was not Lie algebra valued but valued in the
enveloping algebra. Such a condition depends very much on the form
of $A$ and of the calculus $\Omega^2(M)$ and $\Omega^1(H)$. One
could view it as some kind of `differentiability' condition on $A$.
Next we clarify the geometric meaning of our objects.
$\nabla\wedge$ denotes applying the covariant derivative $\nabla$
and then projecting to $\Omega^2(M)$.

\begin{corol} The curvature $R:\Omega^1(M)\to
\Omega^2(M)\tens_M\Omega^1(M)$ for a regular connection obeys  \[
R=((\id\wedge\nabla)-(\extd\tens\id))\circ\nabla.\] The torsion
$T:\Omega^1(M)\to \Omega^2(M)$ and cotorsion $\Gamma\in
\Omega^2(M)\tens_M\Omega^1(M)$ corresponding to
\[ \bar D_A e(v)=\extd e(v)+A(\tilde\pi_{\Omega_0}v\bz)\wedge e(v\bo),
\quad  D_A e^*(w)=\extd e^*(w)+e^*(w\bo)\wedge A(\tilde
\pi_{\Omega_0} w\bt)\]
respectively (assuming a $V$-cobein in the second case)
\[ \nabla\wedge =\extd -T,\quad \Gamma=(\nabla\wedge\id
-\id\wedge\nabla)g+(T\tens\id)g.\]
\end{corol}
\proof These results follow from the general theory outlined in
Section~2 specialised to the bundle $P=M\tens H$ along the lines
already given. However, for trivial bundles one may give a direct
self-contained proof as well. For the curvature the notation
$(\id\wedge\nabla)$ means to act in the second tensor factor of
$\Omega^1(M)\tens_M\Omega^1(M)$ and then project the first two of
the resulting three factors to $\Omega^2(M)$. From the definition
of $\nabla$ we have on a 1-form $\alpha$, \begin{eqnarray*}
R\alpha&&=((\id\wedge\nabla)-(\extd\tens\id))(\extd\alpha^a\tens_M
e(e_a) -\alpha^aA(\tilde\pi_{\Omega_0}e_a\bz)\tens_M e(e_a\bo)\\
&&=\extd \alpha^a\wedge (-A(\tilde\pi_{\Omega_0}e_a\bz)\tens_M
e(e_a\bo)) +A(\tilde\pi_{\Omega_0}e_a\bz)\wedge
A(\tilde\pi_{\Omega_0}e_a\bo\bz)\tens_M e(e_a\bo\bo)\\
&&\quad+\extd (\alpha^a A(\tilde\pi_{\Omega_0}e_a\bz))\tens_M
e(e_a\bo)\\ &&=\alpha^a F_A(\tilde\pi_{\Omega_0} e_a\bz)\tens_M
e(e_a\bo)\end{eqnarray*} using the Leibniz rule and the left
comodule property. This also gives the way to  compute the action
of $R$ from $F_A$. For torsion we project the definition of
$\nabla$ down to $\Omega^2(M)$, so that \begin{eqnarray*}  \nabla
\wedge\alpha&&=(\extd \alpha^a)\wedge e(e_a)-\alpha^a A(\tilde
\pi_{\Omega_0}e_a\bz)\wedge e(e_a\bo)\\ &&=\extd
\alpha-\alpha^a\extd
e(e_a)-\alpha^aA(\tilde\pi_{\Omega_0}e_a\bz)\wedge e(e_a\bo)
=\extd\alpha-\alpha^a\bar D_A e(e_a)\end{eqnarray*}  by the
Leibniz rule in $\Omega^2(M)$. This also makes it clear how $T$
can be  efficiently determined from $\bar D_Ae$. For the cotorsion
we use the metric to similarly relate it to $D_Ae^*$, namely
\begin{eqnarray*}\Gamma&&=D_Ae^*(f^a)\tens_M e(e_a)=\extd
e^*(f^a)\tens_M e(e_a)+e^*(f^a\bo)\wedge
A(\tilde\pi_{\Omega_0}f^a\bt) \tens_M e(e_a)\\ &&=\extd
e^*(f^a)\tens_M e(e_a)+e^*(f^a)\wedge A(\tilde\pi_{\Omega_0}
e_a\bz) \tens_M e(e_a\bo)\\ &&=\extd e^*(f^a)\tens_M
e(e_a)-e^*(f^a)\wedge\nabla\tens_M e(e_a)\end{eqnarray*} where we
use that the right coaction on $V^*$ is adjoint to the left one on
$V$  (obtained as in (\ref{DeltaL})). Specifically, it means that
$f^a\bo\tens e_a\tens f^a\bt=f^a\tens e_a\bo\tens S^{-1}e_a\bt$ for
the relation between the two coactions. Finally, we use the
characterisation of torsion already obtained.
\eproof

The corollary shows in particular one of the key ideas in our
approach\cite{Ma:rie}; the vanishing of cotorsion (or rather the
difference between the torsion and the cotorsion) is a
skew-symmetrized version of the `Levi-Civita' condition of metric
compatibility. From the Riemman tensor above, it is clear that if
we are given a bimodule map
$i:\Omega^2(M)\to\Omega^1(M)\tens_M\Omega^1(M)$ (preferably
splitting the surjection $\wedge$ but not necessarily) we have a
well-defined Ricci tensor \eqn{RicciH}{ {\rm
Ricci}=\<i(R)(e(e_a)),f^a\>=i(F_A(\tilde\pi_{\Omega_0}
e_a\bz))^{ab}e(e_b)\tens_M e(e_a\bo)} where
$F_A=F_A^{ab}e(e_a)\tens_M e(e_b)$ defines its components. The
first trace expression is with the pairing applied to the first
component of $i(R)$ with all coefficients taken to the left in the
$V$-bein basis and $\<me(e_a),f^b\>=m\delta_a{}^b$. It is
independent of the basis of $V$. One may go further and similarly
contract to the scalar curvature. Finally, let us note that we are
taking a view in which the underlying variables are a $V$-bein for
the framing and, given this, an independent $V$-cobein $e^*$ for
the metric. If we fix a specific reference choice of that, e.g.
$e_{ref}^*(f^b)= e(e_a)\eta^{ab}$ for some fixed equivariant
isomorphism $\eta:V^*\isom V$, then any other $V$-cobein has the
form $e^*(f^b)=e^*_{ref}(f^a)g_a{}^b$ for some $g\in GL(n,M)$
where $n=\dim(V)$. Then (summations understood) \eqn{metricM}{
g=e_{ref}^*(f^a)g_a{}^b\tens_M e(e_b)=e(e_a)g^{ab} \tens_M
e(e_b);\quad g^{ac}=\eta^{ab}g{}_b{}^c.}

This completes our treatment of parallelizable quantum Riemannian
manifold structures on general algebras $M$, which can be expected
to be the minimum level of generality for comparison with quantum
theory. The rest of the paper is devoted to constructing examples
of this including quantum groups, finite sets and finite groups.
One could in principle also apply it to specific quantum systems as
well as to discrete algebras such as quaternions in the setting of
\cite{Con:rea}.

\section{Riemannian geometry of quantum groups}

In this section we construct quantum Riemannian geometries where
$M$ is a quantum group. This covers both finite groups and Lie
groups (in an algebraic form) as well as their q-deformations. In
fact Hopf algebras have been used historically to unify Lie theory
and finite group theory and we do the same here by working with
general Hopf algebras. The main result follows in Section~4.1 with
the construction of a natural metrics on the standard $\C_q[G]$
from a braided-Killing form on the braided-Lie algebra tangent to
the fibres of the frame bundle.

For framing we take the same quantum group $H=M$. The classical
meaning of this is explained in \cite{Ma:rie}, with the same
bicovariant differential calculi on $M$ and $H$. These are
determined by ideals $Q_M=Q_H$ as usual. Here
$V=\Omega_0=\ker\eps/Q_H$ has a right coaction $\Ad_R$ and is the
dual of the braided-Lie algebra in the fibre direction. We begin by
checking the various conditions needed to establish a framing or
quantum manifold structure in the sense of Section~3. In effect we
are able for the first time properly to interpret the well-known
`Maurer-Cartan' form in \cite{Wor:dif} in a geometrical manner. It
also provides an actual connection (generally with torsion).

\begin{lemma} For $P=M\tens H$ and $M$ a Hopf algebra, if
$\Omega^1(M)$ is bicovariant then so is $\Omega^1(P)$,
\[ Q_P=Q_M\tens H+M\tens Q_H+\ker\eps_M\tens\ker\eps_H,\quad
\Omega_{0P}=\Omega_{0M}\tens 1\oplus 1\tens\Omega_{0H}\] and the
exterior algebras $\Omega^\cdot(P)$, $\Omega^\cdot(M)$ obey
(\ref{deg2P}). In the case $M=H$ the Maurer-Cartan form
\[ e:\Omega_0\to \Omega^1(H),\quad e(v)
=\pi_{N_H}(S\tilde v\o\tens \tilde v\t)\]
for any representative $\tilde v$ of $v\in\Omega_0$ provides a
framing as well as a zero curvature gauge field
\[ A=e:\Omega_0\to \Omega^1(H).\]
\end{lemma}
\proof For the differential calculus, it is evident that
$\Delta_{M\tens H}(m\tens h)=m\o\tens h\o\tens m
\t\tens h\t$ is a left or right coaction on $M\tens H$ and that
$N_P$ is bicovariant just because $N_M$ and $N_H$ are. The map
$\ver_{M\tens H}$ (not to be confused with that of the bundle)
easily computes as an isomorphism  \[ N_P\isom M\tens Q_M\tens
H\tens H+M\tens M\tens H\tens Q_H+M\tens\ker \eps_M\tens H\tens
\ker\eps_H=M\tens H\tens Q_P\] under the usual identification of
the vector spaces. Note also that we have \[ Q_P=Q_M\tens 1\oplus
1\tens Q_H\oplus\ker\eps_M\tens\ker\eps_H\] as right $\Ad_{M\tens
H}$-comodules. We then apply the Woronowicz construction  for
$\Omega^\cdot(P)$, $\Omega^\cdot(M)$. Here the additional relations
on $\Omega^1(P)\tens_P\Omega^1(P)$ are defined by the kernel of
$\id-\Psi$ where the braiding $\Psi$ is determined by the usual
flip on left and right invariant forms on $P$. But these are just
the images of those either from $M$ or from $H$. Next, that $e$
provides an $\Omega_0$-bein and hence a framing is precisely the
geometric meaning of the isomorphism $\Omega^1(H)\isom
H\tens\Omega_0$, namely with inverse being $s_e$ for the
Maurer-Cartan form. Regularity of $A$ is also immediate since $e$
is known to obey the well-known `Maurer-Cartan equation' \eqn{mc}{
\extd e(v)+e(\tilde\pi_{\Omega_0}\tilde v\o)\wedge
e(\tilde\pi_{\Omega_0}\tilde v\t)=0.} (This in turn is immediate by
working in the universal calculus where $e(\tilde v)=S\tilde
v\o\tens
\tilde v\t-1\tens 1\eps(\tilde v)=S\tilde v\o\extd \tilde v\t$).
From the Maurer-Cartan equation it follows that  if we view $A=e$
as a gauge field then it is regular and has zero curvature.
\eproof

The operators $\rho_a{}^b$ in (\ref{partialM}) for this framing
are those of right translation according to \eqn{emod}{ e(v)g=g\o
e(vg\t),\quad\forall v\in\Omega_0,\quad g\in H.} There is also a
right-handed framing defined by $\bar e(v)=\pi_{N_H}(\tilde
v\o\tens S\tilde v\t)$ and related by \eqn{eebar}{ e(v)=\bar
e(S\pi_{\Omega_0}\tilde v\t)(S\tilde v\o)\tilde v\th.} Hence the
braiding $\Psi$ in the definition $\Omega^\cdot(H)$ can be written
in the crossed module form \eqn{PsiD}{ \Psi(e(v)\tens_H
e(w))=e(\pi_{\Omega_0}\tilde w\t)\tens_H e(v(S\tilde w\o)\tilde
w\th)} rather than the more standard form with $e,\bar e$ as in
(\ref{worskew}). We clearly have a natural `lift' \eqn{iwor}{
i=\id-\Psi:\Omega^2(H)\to \Omega^1(H)\tens_H\Omega^1(H),} since
$\Omega^2(H)$  is by definition $\Omega^1(H)\tens_H\Omega^1(H)$
modulo $\ker(\id-\Psi)$ and hence isomorphic to the image of
$\id-\Psi$. On the other hand, $\Psi$ does not generally obey
$\Psi^2=\id$ and as a result this map does {\em not} generally
split $\wedge$, i.e. $i\circ\wedge$ is not a projection. Therefore
one can use this $i$ to define the Ricci tensor and interior
products, etc., but it is not necessarily the best choice.

The torsion tensor corresponds from Section~3 to
 \eqn{tore}{ D_Ae(v)=\extd
e(v)+ A(\tilde\pi_{\Omega_0}(S^{-1}\tilde v\th)\tilde v\o) \wedge
e(\pi_{\Omega_0}\tilde v\t)} since the coaction on $\Omega_0$ to be
used is the right adjoint one converted to a left coaction by
(\ref{DeltaL}). We do not solve this in general (this would appear
to require further data) but it is worth noting that classically
$A=\h e$ is a torsion free connection, and also cotorsion free for
the Killing metric for any classical compact Lie
group\cite{Ma:rie}. The latter is an example of an important class
of quantum metrics where $\Omega_0$-cobein $e^*:\Omega_0^*\to
\Omega^1(H)$ is defined by a nondegenerate $\Ad$-invariant bilinear form.
Such an element corresponds to an $\Ad$-invariant element of
$\eta=\eta\uo\tens\eta\ut \in\Omega_0\tens\Omega_0$ nondegenerate
as a map $\eta:\Omega_0^*\to\Omega_0$ by evaluation against the
second component.

\begin{propos} Any nondegenerate $\Ad$-invariant $\eta\in\Omega_0
\tens\Omega_0$ defines a coframing $e^*=e\circ\eta$. The
corresponding metric $g=e(\eta\uo)\tens_H e(\eta\ut)$ is symmetric
in the sense $\wedge(g)=0$ {\em iff} $\eta=\eta\ut\tens
S^2\eta\uo$. Its cotorsion in terms of $e$ is given by
\[ D_Ae(v)=\extd e(v)+ e(\pi_{\Omega_0}\tilde v\t)\wedge
A(\tilde\pi_{\Omega_0}(S\tilde v\o)\tilde v\th).\] \end{propos}
\proof For the framing the only delicate part is to check
that $\eta:\Omega_0^*\to\Omega_0$ is equivariant, where
$\Omega_0^*$ has the right coaction adjoint to the left coaction on
$\Omega_0$ given as in (\ref{DeltaL}) by $S^{-1}$, i.e. that
$\eta\uo\tens\eta\ut\t\tens S^{-1}((S\eta\ut\o)\eta\ut\th)
=\eta\uo\t\tens\eta\ut\tens (S\eta\uo\o)\eta\uo\th$, using Hopf algebra
methods\cite{Ma:book}. Next, the condition that $e(\eta\uo)\wedge
e(\eta\ut)=0$ is that $(e\tens_H e)\circ\eta$ is in the kernel of
$(\id-\Psi)$ where $\Psi$ is as above. The corresponding $\Psi$ on
$\Omega_0\tens\Omega_0$ computes as
\[\Psi(\eta\uo\tens\eta\ut)=\eta\ut\t\tens
\eta\uo(S\eta\ut\o)\eta\ut\th =\eta\ut\tens S^2 \eta\uo\]
in view of the equivariance of $\eta$. Finally, cotorsion from
Section~3 corresponds to
\[ D_A e^*=(\extd\tens\id) e^*+e^*\uo \wedge
A(\tilde\pi_{\Omega_0}e^*\ut\bz) \tens e^*\ut\bo
\] where $e^*=e^*\uo\tens e^*\ut\in
\Omega^1(H)\tens\Omega_0$, or equivalently \eqn{cotorH}{
D_Ae^*(w)=\extd e^*(w)+e^*(w\bo)\wedge A(\tilde\pi_{\Omega_0}w\bt)}
for the right coaction on $\Omega^*_0$ adjoint to the left coaction
on $\Omega_0$ obtained as in (\ref{DeltaL}). This second form where
$e^*:\Omega^*_0\to\Omega^1(H)$ and equivariance of $\eta$
immediately gives the result.
\eproof

Hence we have a canonical framing and metric and at least one
natural (not generally torsion free or cotorsion free) connection
on any Hopf algebra, and concrete equations for the torsion and
cotorsion conditions. We also have a `tautological' choice of
`gamma' matrix and hence an induced Dirac operator for each
connection. Thus, let $W$ be a right $H$-comodule viewed as in
(\ref{DeltaL}) as a left comodule. Also let the inverse map
$\eta^{-1}(v)=\eta\umo\<\eta\umt,v\>$ define
$\eta^{-1}\in\Omega_0^*\tens\Omega_0^*$ or
$\eta^{-1}:\Omega_0\tens\Omega_0\to k$ depending on ones point of
view (we assume finite-dimensionality).

\begin{corol}  For any right comodule $W$ and
$\eta$ as above there is a canonical equivariant map
\[ \gamma:\Omega_0\tens W\to W,\quad \gamma(v)w
=\eta^{-1}(v\tens \tilde\pi_{\Omega_0}(w\bz))w\bo \]
obeying additionally the identity
\[ (\gamma\circ\gamma)(\eta)w=\<c,w\bz\>w\bo,\quad c=\eta\umo\eta\umt.\]
\end{corol} \proof  By similar Hopf algebra methods equivariance of
$\eta$ can be written as $\eta^{-1}(v\bo\tens w\bo)v\bt
w\bt=\eta^{-1}(v\tens w)$. From this one similarly computes
\[ \Delta_R(\gamma(v)w)
=w\bo\bo\eta^{-1}(v\bo\tens S^{-1}w\bo\bt)\tens v\bt
w\bt=\gamma(v\bo)w\bo\tens v\bt w\bt,\]
\begin{eqnarray*}
\gamma(\eta\uo)\gamma(\eta\ut)w&&= w\bo\bo\eta^{-1}(\eta\uo\tens
S^{-1}w\bo\bt)\eta^{-1}(\eta\ut\tens S^{-1}w\bt)\\
&&=w\bo\bo\eta^{-1}(S^{-1}w\bt\tens
S^{-1}w\bo\bt)=w\bo\eta^{-1}((S^{-1}w\bt)\o\tens (S^{-1}w\bt)\t)
\end{eqnarray*} as required. Note that $c$ is invariant
under the right coadjoint coaction on $\Omega_0^*$ because
$\eta\umt\tens\eta\umo$ is (the reversal is because it is the left
coadjoint coaction that respects the product here).   \eproof

There is also a tautological $\gamma^*$ defined similarly without
the $\eta^{-1}$ i.e. just from the comodule itself and with similar
features. The equivariance of $\gamma$ (and $\gamma^*$) here
replaces the idea that the antisymmetric products of  $\gamma$
classically generates a representation of the rotation group or
that $\gamma$ generates a representation of the spin group.
Meanwhile, the coadjoint invariant element $c$ is central at least
when it lies in a Hopf algebra $U$ dually paired with $H$ (which
will generally be the case). We denote by $\rho_W$ the left action
of $U$ corresponding to the right coaction of $H$ so when $\rho_W$
is irreducible then $(\gamma\circ\gamma)(\eta)$ etc. will be a
multiple of the identity, which is a remnant of the usual `Clifford
algebra' property for the symmetric products of $\gamma$.

\begin{propos} With framing and connection provided by the Maurer-Cartan
form itself and with the tautological $\gamma$ as above, the Dirac
operator associated to any right $H$-comodule is
\[ \Dsl=\del^a \gamma_a-\rho_W(S^{-1}c),\quad
\gamma_a=\eta^{-1}_{ab}\rho_W(S^{-1}f^b).\]   \end{propos} \proof
Here $\eta^{-1}_{ab}=\eta^{-1}(e_a\tens e_b)$. The general
expression for the Dirac operator  is in Corollary~3.6. We note
that if $A=e=A^ae(e_a)$ then its components are
$A^a(v)=\<f^a,v\>$. Here $f^a$ are a dual basis of
$\Omega_0^*\subset\ker\eps\subset U$ (which we assume for
convenience of presentation). Hence
$\<f^a,h\>=\<f^a,\tilde\pi_{\Omega_0}h\>$ automatically makes the
projection, giving the general form of $\Dsl$ as stated. We write
the coaction as an action of the dual basis for convenience. For
the particular form of $\gamma$ itself given by the coaction or by
$\rho_W$ we immediately obtain the result stated. \eproof

This completes our analysis for general Hopf algebras. Before
turning to nontrivial examples let us note that for $H$
cocommutative (e.g. classically an Abelian group) all connections
$A$ are torsion free and induce the same $\nabla$ given by
$\nabla\alpha=\extd\alpha^a\tens_H e(e_a)$ with zero Riemannian
curvature. Any nondegenerate bilinear form
$\eta\in\Omega_0\tens\Omega_0$ defines a metric with zero
cotorsion as well. This does however, give a simple example of
noncommutative geometry fully in keeping with the classical
picture. For example, for a Lie algebra $\cg$ the enveloping
algebra $H=U(\cg)$ can be viewed `up side down' as the
quantisation of the Kirillov-Kostant bracket on $\cg^*$.

\begin{propos} For $H=U(\cg)$, coirreducible calculi are provided by
$(V,\lambda)$ with $V$ an irreducible right module (with right
action $\rho$) and $\lambda\in P(V)$ a ray. Here
\[ \Omega_0=\ker\eps/\ker\rho_\lambda,\quad
\rho_\lambda:\Omega_0\isom V,\quad \rho_\lambda(h)=\lambda\cdot\rho(h),
\quad\forall h\in U(\cg).\]
Then $e=e_{MC}\circ\rho_\lambda^{-1}$ is a framing, where $e_{MC}$
is the  Maurer-Cartan form, and
\[ e(v)\xi=\xi e(v)+e(v\cdot\rho(\xi)),\quad
\extd(\xi_1\cdots\xi_n)=\sum_{m=0}^{n-1}\sum_{\sigma\in
S_{m|n-m}}\xi_{\sigma(1)}\cdots\xi_{\sigma(m)}
e(\rho_\lambda(\xi_{\sigma(m+1)}\cdots\xi_{\sigma(n)})),\] where
$\xi,\xi_i\in \cg$ and $S_{m|n-m}$ denotes permutations of
$\{1,2,\cdots,n\}$ such that $\sigma(1)<\cdots<\sigma(m)$ and
$\sigma(m+1)<\cdots<\sigma(n)$ (an $m$-shuffle). Any bilinear form
$\eta$ in $V\tens V$ defines a metric as above, and $\nabla$ is
torsion free and cotorsion free.
\end{propos} \proof The differential
calculus is a `differentiation' of the classification in
\cite{Ma:cla} for the calculi for group algebras as a pair
consisting of an irreducible representation and ray.  After
differentiating those formulae one verifies directly that the above
defines a calculus and that it is coirreducible. Here
$\ker\rho_\lambda$ is clearly an ideal and for fixed $\rho$ and in
the irreducible case the image of $\rho_\lambda$ must be all of
$V$. Actually the minimum we need for a calculus here is that
$\lambda$ is a cyclic vector. If we simply identify $\Omega_0$ with
$V$ in this way then clearly \eqn{extdlie}{\extd
\xi=\lambda\rho(\xi),\quad v\xi=\xi v+v\rho(\xi)} which is easily
seen to extend by Leibniz to a well-defined calculus. Thus
\[ \extd
(\xi\eta)=(\lambda\rho(\xi))\eta+\xi(\lambda\rho(\eta))
=\xi\lambda\rho(\eta)+\eta\lambda\rho(\xi)+\lambda\rho(\xi\eta)\]
so that $\extd(\xi\eta-\eta\xi)=\extd[\xi,\eta]$. A proof by
induction gives the general form of $\extd$  (writing the
identification $e$ explicitly). Also the right action on $V$
corresponds to right multiplication on $\Omega_0$ as it should,
since
$\rho_\lambda(h\xi)=\lambda\rho(h\xi)=\lambda\rho(h)\rho(\xi)
=\rho_\lambda(h)\rho(\xi)$. If $\Omega_0'$ defines a quotient
differential calculus then it corresponds to a surjection
$\phi:V\to \Omega_0'$ an intertwiner as $U(\cg)$-modules which,
for irreducible $V$, must be an isomorphism. To form a commutative
triangle,
$\extd\xi=\phi(\lambda\rho(\xi))=\phi(\lambda)\rho'(\xi)$, say, so
that the quotient calculus is isomorphic to our $(V,\lambda)$
calculus with $\lambda'=\phi(\lambda)$. Moreover, $(V,\lambda)$ is
isomorphic to $(V,\lambda')$ if and only if $\phi$ is a nonzero
multiple of the identity i.e. $\lambda'$ proportional to
$\lambda$, i.e. the calculus depends on $\lambda$ only up to
scale. This describes the calculus that we use. While these are
not all possible calculi (any ideal in $\ker\eps$ defines a
calculus since $H$ is cocommutative), they are the natural
`integrable' calculi in the sense that they `differentiate' the
formulae in the finite group case. We compute the geometric
structure. This is defined in terms of $\Omega_0$ (which is hard
to work with) so we work instead with the its isomorphic image
which is $V$ as stated. Hence we take $V$ itself as the framing
space and $e$ the Maurer-Cartan form converted under the
identification (similarly for all the formulae above). For the
exterior algebra we have $\extd e(v)=0$ and $e(v)\wedge
e(w)=-e(w)\wedge e(v)$. \eproof

More generally, it is clear from the proof that any representation
$V$ and cyclic $\lambda$ likewise gives a framing, etc. (if we do
not care about irreducibility). This describes $U(\cg)$ as a
`noncommutative flat space' (namely quantized $\cg^*$). One can
also choose interesting spinor spaces and $\gamma$-matrices and
hence a Dirac operator sensitive to $A$. On $U(su_2)$ for example
one could take the usual $\gamma$ (Pauli) matrices. And, of course,
one can have other metrics not induced by constant $\eta$.

\subsection{Killing form metric on $\C_q[G]$}

We now turn to our main construction which is the example of $M=H$
a dual quasitriangular Hopf algebra. It means that there is a
`universal R-matrix functional' $\CR:H\tens H\to k$, which
includes the standard deformations $\C_q[G]$ of the classical
simple Lie groups. $\Omega^1(H)$ is built from a
finite-dimensional right comodule $W$ (which we view as a left
module of $H^*$ with action $\rho_W$. The element
$\CQ=\CR_{21}\CR$ is the `universal Killing form' and we view it
as a map  $\CQ:H\to H^*$ by evaluation, i.e.
$\<g,\CQ(h)\>=\CQ(h\tens g)=\CR(g\o\tens h\o)\CR(h\t\tens g\t)$
for $g,h\in H$. We assume that $\rho_W\circ\CQ$ is surjective
(e.g. if $\CR$ is factorisable and $\rho_W$ irreducible). We also
define the induced actions of $H$: \eqn{rhopm}{
\rho_+(h)^\alpha{}_\beta=\CR(h\tens\rho^\alpha{}_\beta),\quad
\rho_-(h)^\alpha{}_\beta=\CR^{-1}(\rho^\alpha{}_\beta\tens h)}
where $e_\alpha\mapsto e_\beta\tens \rho^\beta{}_\alpha$ defines
the matrix elements   of $\rho_W$ for a basis $\{e_\alpha\}$ of
$W$. With these notations one knows that there is a bicovariant
differential calculus defined by \eqn{OmegaGq}{
\Omega_0=\ker\eps/\ker\rho_W\circ \CQ,\quad
\rho_W\circ\CQ:\Omega_0\isom \End(W).} This  is part of the
construction in \cite{Ma:cla}, where it was shown that such
calculi with $\rho_W$ irreducible essentially classify all the
coirreducible calculi for factorisable quantum groups such as
$\C_q[G]$. We let $W^\circ$ be the predual of $W$ as a right
comodule.

\begin{propos} A dual-quasitriangular Hopf algebra $H$ with
calculus defined by $(W,\rho_W)$ is framed by $V=\End(W)=W\tens
W^\circ$ and $e=e_{MC}\circ(\rho_W\circ \CQ)^{-1}$. We have
\[ e(\phi)h=h\o e(\rho_-(Sh\t)\circ\phi\circ\rho_+(h\th)),\quad \extd
h=\cdot(\id\tens e)(h\o\tens \rho_W\circ\CQ(h\t)-h\tens\id)\] for
all $h\in H$, $\phi\in \End(W)$. Moreover, there is a natural
choice of spinor space, namely $W$, with equivariant $\gamma:V\tens
W\to W$ provided by the identity matrix and
\[(\Dsl\psi)^\alpha=\del^\alpha{}_\beta\psi^\beta-
A(\tilde\pi_{\Omega_0}S^{-1}
\rho^\beta{}_\gamma)^\alpha{}_\beta\psi^\gamma,\] where
$\psi^\alpha\in H$ are the spinor components.
\end{propos}  \proof  With the identification
(\ref{OmegaGq}) understood, one could write the calculus
$\Omega^1(H)$ as \eqn{extdQ}{\extd h= h\o\rho_W\circ
\CQ(h\t)-h\id,\quad \phi h=h\o
\rho_-(Sh\t)\circ\phi\circ\rho_+(h\th),\quad\forall h\in H} for
the exterior derivative and bimodule structure on
$\phi\in\End(W)$. In our context this identification is made by
the framing and gives the structure shown when we write this
explicitly. One may check that
 \[
\rho_W\circ\CQ(hg)=\rho_-(S(hg)\o)\rho_+((hg)\t)=\rho_-(Sg\o)\rho_-(Sh\o)
\rho_+(h\t)\rho_+(g\t)
 =\rho_-(Sg\o)\rho_W\circ\CQ(h)\rho_+(g\t)\] which leads to the stated
$H$-module structure on $V$. Meanwhile, the right adjoint coaction
is known\cite{Ma:book} to intertwine under $\CQ$ with the right
coadjoint coaction on $H^*$, which means
\eqn{rhoQAd}{\rho_W\circ\CQ(h\t)^\alpha{}_\beta\tens
(Sh\o)h\th=\rho_W\circ\CQ(h)^a{}_b\tens \rho^\alpha{}_a
S\rho^b{}_\beta.} In our present setting the equivariance follows
easily from the dual-quasitriangularity axioms for $\CR$ provided
the coaction maps a dual basis element as $f^\alpha\mapsto f^\beta
\tens S\rho^\alpha{}_\beta$. This means that we identify $V=W\tens
W^\circ$ as stated. It is straightforward to verify that $\extd$
as stated obeys the Leibniz rule and that we indeed have a
calculus. Also from (\ref{rhoQAd}) and (\ref{PsiD}) one obtains
easily the braiding in terms of R-matrices
$R^\alpha{}_\beta{}^\gamma{}_\delta=\CR(\rho^\alpha{}_\beta\tens
\rho^\gamma{}_\delta)$ and $\tilde
R^\alpha{}_\beta{}^\gamma{}_\delta=\CR(\rho^\alpha{}_\beta\tens
S\rho^\gamma{}_\delta)$, \eqn{PsiGq}
{\Psi(\phi\tens\psi)^\alpha{}_\beta{}^\gamma{}_\delta
=R^a{}_{\mu}{}^\alpha{}_b\phi^\mu{}_\nu
R^b{}_\sigma{}^\nu{}_c\psi^\sigma{}_\tau
R^{-1}{}^\tau{}_d{}^c{}_\delta\tilde
R^\gamma{}_a{}^d{}_\beta,\quad \forall \phi,\psi\in \End(W),} or
$\sum_i \psi^i_2R\phi^i_1 R_{21}=R\phi_1R_{21}\psi_2$ if
$\Psi(\phi\tens\psi)\equiv \sum_i \psi^i\tens\phi^i$ in a standard
notation. The partial derivatives are as usual the coefficient in
$\extd$ of the basic forms $e(e_\alpha\tens f^\beta)$, which means
\eqn{delQ}{ \del^\alpha{}_\beta
(h)=h\o\rho_W\circ\CQ(h\t)^\alpha{}_\beta-h\delta^\alpha{}_\beta.}
We define the gamma-matrices as projectors in $W\tens W^\circ$
acting by evaluation (or $\gamma:V\to W\tens W^\circ$ the identity
map). Thus $\gamma_\beta{}^\alpha(\psi)=e_\beta \psi^\alpha$ and
$\gamma_\beta{}^\alpha=e_\beta\tens f^\alpha$, giving $\Dsl$ as
stated.  \eproof

\begin{figure}
\[\epsfbox{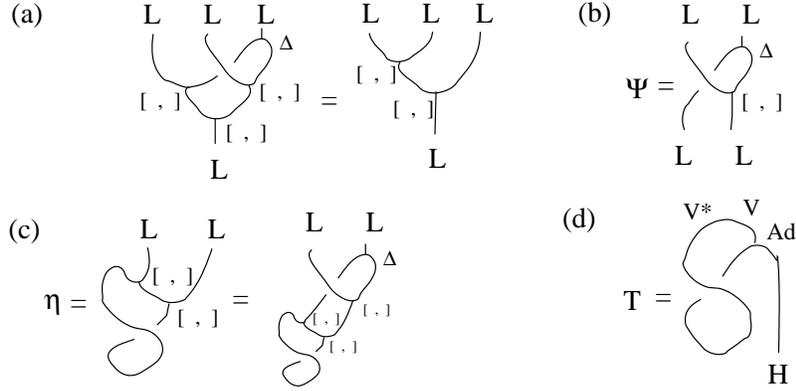}\]
\caption{Pentagonal axiom (a) of a braided-Lie algebra. Induced
`double braiding' (b), braided-Killing form and its braided
symmetry (c) and construction (d) as $\eta=\Delta T$ in our
formulation.}
\end{figure}

We turn now to the construction of a natural `Killing form'
metric. In fact we will give a self-contained quantum-group
construction which avoids braided categories, but the following is
the picture behind it. Thus, it was shown in \cite{Ma:cla} that
for all such calculi the dual of $\Omega_0$ forms a braided-Lie
algebra $L$ in the sense of \cite{Ma:lie}. These are modelled on
the properties of the 1-dimensional extension $\cg\oplus k.c$ when
$\cg$ is an ordinary Lie algebra; there is a coproduct
$\Delta:L\to L\tens L$ and an extended bracket $[\ ,\ ]:L\tens
L\to L$ and everything lives in a braided category (classically we
would extend by $[\xi,c]=\xi$, $[c,\xi]=0$ for $\xi\in\cg$ and
$[c,c]=c$ with $\Delta c=c\tens c$, $\Delta\xi=\xi\tens
c+c\tens\xi$ for the coproduct, and have a trivial braiding). The
main `pentagonal Jacobi identity' axiom of a (right-handed)
braided Lie algebra is shown in Figure~1(a) in a diagrammatic
notation\cite{Ma:book} with operations flowing down the page and
with the braid-crossing denoting the `background braiding' of the
category. The axioms for $(L, [\ ,\ ],\Delta)$ and a counit are
strong enough to define an additional `double' braiding $\Psi$
shown in Figure~1(b) and from this an enveloping algebra $U(L)$ as
a bialgebra or `braided group' in the braided category. This is
defined as the quadratic algebra generated by $L$ with relations
of symmetry with respect to $\Psi$ (i.e. setting to zero the image
of $\id-\Psi$) and coproduct extending $\Delta$ on $L$
(classically this would recover a quadratic extension of the usual
$U(\cg)$). There is also a braided-Killing form $\eta$ in
Figure~1(c) which is shown there to be braided-symmetric in the
sense $\eta=\eta\circ\Psi$. Here $\cup$ and $\cap$ are evaluation
and coevaluation of $L$ with a suitable dual. The braided-Killing
form $\eta$ classically restricts to the usual one on $\cg$ and
$\eta(c,c)=1$. Thus Lie theory is contained as a special class of
braided-Lie algebras and acquires extra structure such as the
double braiding $\Psi$.

In our case $L=W^*\tens W=V^*$ in the preceding proposition with
basis $\{x^\alpha{}_\beta=f^\alpha\tens e_\beta\}$ and $\Delta
x^\alpha{}_\beta=x^\alpha{}_\gamma\tens x^\gamma{}_\beta$ has a
matrix form. The Lie bracket $[\ ,\ ]$ is given in \cite{Ma:lie}
in an R-matrix form as well as the background braiding defined by
$\CR$. The double braiding $\Psi$ is the adjoint of (\ref{PsiGq})
for the exterior algebra and correspondingly the enveloping
bialgebra $U(L)$ is the left-handed braided matrices $B_L(R)$ with
relations $ x_2Rx_1 R_{21}=Rx_1R_{21}x_2$. There is an algebra map
$U(L)\to H^*$ sending $x^\alpha{}_\beta$ to $\rho_W\circ\CQ(\
)^\alpha{}_\beta\in H^*$ and we identify the image of $L$ with
$\Omega_0^*$ by the counit projection to
\eqn{fab}{f^\alpha{}_\beta=\rho_W\circ\CQ(\
)^\alpha{}_\beta-\delta^\alpha{}_\beta\eps\in H^*} adjoint to the
restriction to $\ker\eps\subset H$ in (\ref{OmegaGq}). The braided
Killing form on $L\tens L$ can then be viewed in
$\Omega_0\tens\Omega_0$. We now give a version of this
construction directly in our setting.

\begin{theorem} Let $H$ be a dual-quasitriangular Hopf algebra with
differential calculus as above. There is a braided-symmetric and
$\Ad$-invariant `braided-Killing form'
\[\eta_{\Omega_0}=
(\tilde\pi_{\Omega_0}\tens\tilde\pi_{\Omega_0})\Delta T\in
\Omega_0\tens\Omega_0;\quad T=\CR(\tau\o\tens S\tau\t)
\tau\th,\quad \tau=\rho^\alpha{}_\alpha S\rho^\beta{}_\beta\in H\]
\[ \eta=(\rho_W\circ \CQ\tens\rho_W\circ\CQ)(\Delta T)
-\id\tens\rho_W\circ\CQ(T)-\rho_W\circ\CQ(T)\tens\id
+\eps(T)\id\tens\id\in V\tens V.\] If nondegenerate, there is a
braided-symmetric Riemannian metric $g=(e\tens_H e)\eta$ with
$\wedge(g)=0$.
\end{theorem}
\proof The two applications  of $[\ ,\ ]$ in the `figure of eight'
braided trace in Figure~1(c) can be written as a product in $U(L)$
followed by a single $[\ ,\ ]$ and this dualises to the coproduct
of $U(L)^*$ applied to the element $T$ in Figure~1(d) (after some
convention adjustments). This coproduct of $U(L)^*$ is essentially
that of $H$, so we have \eqn{etaT}{\eta'=(\rho_W\circ\CQ\tens
\rho_W\circ\CQ)\Delta T\in V\tens V,} where we use $\rho_W\circ Q$
to map $H$ to $V$. This is the natural object from the braided-Lie
theory and we will see that it has the stated features of $\eta$,
however for our geometrical application we have to first project
$\Delta T$ down to $\Omega_0\tens\Omega_0$ which is
$\eta_{\Omega_0}$ as stated (we have done the same in previous
sections in the expressions $(\id\tens
\tilde\pi_{\Omega_0})\Ad:\Omega_0\to \Omega_0\tens\Omega_0$ dual
to $[\ ,\ ]$). Or by (\ref{OmegaGq}) we apply the counit
projection $\pi_\eps$ to $\Delta T$  and then $\rho_W\circ\CQ$ to
give the corresponding element of $V\tens V$.

We now directly verify the properties of $\eta'$ and hence $\eta$.
Notice first that if $\{e_a\}$ is a basis of $V$ and $\{f^a\}$ a
dual basis of $V^*$, let \eqn{tau}{ \tau=\<e_a\bo,f^a\>e_a\bt.}
Cyclicity of the trace here appears as the following fact:
\eqn{cycl}{\tau\o\tens \tau\t \tens\dots\tens
\tau_{\scriptstyle(n)}=\tau_{\scriptstyle(n)}\tens
\tau\o\tens\cdots\tens \tau_{\scriptstyle(n-1)}.} This is because
the first expression may be written as the trace of $n$
applications of the coaction,
\[ \tau=\<e_a\bo{}^{\cdots}\bo,f^a\>e_a\bo{}^{\cdots}\bo\bt\tens
\cdots\tens e_a\bo\bt\tens e_a\bt.\] The outermost $\bo$ is
equivalent due to equivariance of the duality pairing to $f^a\bo$
and $e_a\bo{}^{\cdots}\bo\bt$ replaced by $Sf^a\bt$. On the other
hand $f^a\bo\tens e_a\tens Sf^a\bt=f^a\tens e_a\bo\tens e_a\bt$ by
a change of basis (or invariance of the coevalution element
$f^a\tens e_a$). Hence we may replace the coaction on $f^a$ by an
innermost coaction $e_a\bo$, putting an extra $\bo$ in all the
other places and replacing $Sf^a\bt$ by $e\bt$. Converting the
iterated coactions back to coproducts gives the cyclicity property
(in fact one needs only a coalgebra for the cyclicity with the
appropriate adjoint operation in the role of $S$). In our case
$V=W\tens W^\circ$ with the coaction given in the preceding
proposition. Then \[ \tau=\<(e_\alpha\tens
f^\beta)\bo,f^\alpha\tens e_\beta\>(e_\alpha\tens f^\beta)\bt
=\<f^\alpha\tens e_\beta,e_a\tens f^b\>\rho^a{}_\alpha
S\rho^\beta{}_b=\rho^\alpha{}_\alpha S \rho^\beta{}_\beta.\]

Now we compute the figure-of eight braided trace, which is fairly
routine\cite{Ma:book}. We read Figure~1(d) from the top down,
starting with $f^a\tens e_a$. This becomes $f^a\tens e_a\bo\tens
e_a\bt$. We then apply the background braiding to the first two
places and evaluation, to find
\begin{eqnarray*} T&&=\<e_a\bo\bo,f^a\bo\>\CR(f^a\bt\tens e_a\bo\bt)
 e_a\bt=\<e_a\bo\bo\bo,f^a\>\CR(S^{-1}e_a\t\tens
e_a\bo\bo\bt)e_a\bo\bt\\
&&=\tau\t\CR(S^{-1}\tau\th\tens \tau\o)=\CR(\tau\o\tens
S\tau\t)\tau\th=\cv(\tau\o)\tau\t
\end{eqnarray*}
using that $\CR$ is $S\tens S$-invariant and cyclicity
(\ref{cycl}) again. Note that $T=\cv(\tau\o)\tau\t$ where
$\cv(h)=\CR(h\o\tens Sh\t)$ implements $S^{-2}$ by
convolution\cite{Ma:book}. Next,
\begin{eqnarray*} \Ad T&&=\CR(\tau\o\tens S\tau\t)\tau\fo\tens
(S\tau\th)\tau\fiv=\CR(\tau\t\tens S\tau\th)\tau\fiv\tens
(S\tau\fo)\tau\o\\
&&=\CR(\tau\o\tens S\tau\fo)\tau\fiv\tens \tau\t S\tau\th=T\tens 1
\end{eqnarray*}
by cyclicity and the axioms of a dual-quasitriangular structure or
that $\cv$ implements $S^{-2}$\cite{Ma:book}. So $T$ and hence
$\eta'$ are $\Ad$-invariant (since $\Delta$ and $\rho_W\circ\CQ$
(and $\tilde\pi_{\Omega_0}$) are $\Ad$-covariant). Similarly from
the cyclicity (\ref{cycl}) and the property of $\cv$ it is clear
that $ST=T$ and  $(\id\tens S^2)\Delta^{\rm op}T=\Delta T$ so
$\eta'$ and hence $\eta$ are $\Psi$-invariant as in
Proposition~4.2. One also has explicit formulae using the
definition of $\CQ$ and the dual-quasitriangularity axioms for
$\CR$,
\eqn{etamat}{\eta'{}^\alpha{}_\beta{}^\gamma{}_\delta=u^{b_5}{}_{b_6}
Q^{b_6}{}_{b_1}{}^{a_2}{}_{a_3} \tilde
R{}^{a_1}{}_{a_2}{}^{a_3}{}_{a_4} Q^{a_4}{}_{a_5}{}^a{}_b\tilde
R^\alpha{}_a{}^{b_4}{}_{b_5} R^{-1}{}^{b_3}{}_{b_4}{}^b{}_\beta
Q^{a_5}{}_{a_1}{}^c{}_d\tilde R{}^\gamma{}_c{}^{b_2}{}_{b_3}
R^{-1}{}^{b_1}{}_{b_2}{}^d{}_\delta}
\eqn{Tmat}{\rho_W\circ\CQ(T)^\alpha{}_\beta=\tilde
R{}^{a_1}{}_{a_2}{}^{a_3}{}_{a_4}
Q^{b_4}{}_{b_1}{}^{a_2}{}_{a3}Q^{a_4}{}_{a_1}{}^a{}_b
R^{-1}{}^{b_1}{}_{b_2}{}^b{}_\beta \tilde
R{}^\alpha{}_a{}^{b_2}{}_{b_3} u^{b_3}{}_{b_4},\quad
\eps(T)=u^{b_1}{}_{b_2}Q^{b_2}{}_{b_1}{}^{a_2}{}_{a_3}\tilde
R^{a_1}{}_{a_2}{}^{a_3}{}_{a_1},} where $u=\tilde
R{}^a{}_\cdot{}^\cdot{}_a$ and $Q=R_{21}R$. \eproof

The braided-Killing form the standard quantum groups $\C_q[G]$ is
closely related to the usual Killing form and is typically
nondegenerate for generic $q\ne 1$ (being rational functions in
$q$ they need to be nondegenerate at only one point to establish
this). Hence the theorem above provides a construction of the
metric for such quantum groups and their standard bicovariant
differential calculi. We will demonstrate this explicitly for the
case of $\C_q[SU_2]$ with its standard 4-dimensional calculus.
(here $W$ is the spin $\h$ representation). The exterior algebra
in this case is well-known and in our conventions is as follows.
We let $e_1{}^1=e_a$, $e_2{}^1=e_b$, etc. for brevity, and
$\theta=e_a+e_d$. Then $e_a,e_b,e_c$  behave like usual forms or
Grassmann variables and
\[ e_a\wedge e_d+e_d\wedge e_a+\lambda e_b\wedge e_c=0,\quad
e_d\wedge e_b+q^2e_b\wedge e_d+\lambda e_a\wedge e_b=0,\quad
e_c\wedge e_d+q^2e_d\wedge e_c+\lambda e_c\wedge e_a=0,\]
\[e_d^2=\lambda e_b\wedge e_c,\quad \extd =[\theta,\ ],\quad e_a
\begin{pmatrix}a&b\\ c&d\end{pmatrix}=\begin{pmatrix}qa&q^{-1}b\\
qc&q^{-1}d\end{pmatrix}e_a,\quad [e_b, b]=[e_b,
d]=[e_c,a]=[e_c,c]=0\]
\[[ e_b, a]=q\lambda\, be_a,\quad [e_b, c]=q\lambda\, de_a,\quad
[e_c, b]=q\lambda\, ae_a,\quad [e_c, d]=q\lambda\, ce_a\]
\[ [e_d,a]_{q^{-1}}=\lambda be_c,\quad
[e_d,b]_q=\lambda ae_b+q\lambda^2be_a,\quad
[e_d,c]_{q^{-1}}=\lambda de_c,\quad [e_d,d]_q=\lambda
ce_b+q\lambda^2d e_a\] where $[x,y]_q=xy-qyx$ and
$\lambda=(1-q^{-2})$, and $a,b,c,d\in SU_q(2)$. Note that the
$\wedge$ relations are essentially those for the exact
differentials on q-Minkowski space\cite[Sec. 10.5]{Ma:book} given
by the braid statistics $\Psi_+$ for the addition law on that, as
must be the case because $B_L(R)$ is the coordinate algebra of
$q$-Minkowski space as well as $U(L)$ (details will appear
elsewhere \cite{GomMa:coh}).

\note{ e_\alpha{}^\beta t^i{}_j=t^i{}_k
R^\gamma{}_\alpha{}^k{}_l{} e_\gamma{}^\delta
R^l{}_j{}^\beta{}_\delta\] where $t^1{}_1=a$, $t^1{}_2=b$ etc.,
are the usual $\C_q[SU_2]$ generators. The $R$-matrix relations
can easily be computed explicitly as}

\begin{propos} Let
$[n]_q=(1-q^n)/(1-q)$. The braided-Killing form for the spin $\h$
differential calculus on $\C_q[SU_2]$, divided by $(q-1)^2$ is
\[\eta=q^{-12}[8]_q[2]_q\, \eta_K-\lambda
\left([3]_q(q^{-9}+q^{-7})-[2]_q^2q^{-2}\right)\theta\tens
\theta\]
\[ \eta_K=e_b \tens e_c + q^2
e_c\tens e_b + {(e_a \tens e_a - q e_a \tens e_d - q e_d \tens e_a
+ q(q^2 + q - 1)e_d \tens e_d)\over [2]_q}\]
\end{propos}
\proof We use the R-matrix formulae obtained in Theorem~4.7. In
fact the difference between $\eta$ and $\eta'$ is a multiple of
$\theta\tens \theta=\id\tens\id$ so only affects the second term
here. One has $\rho_W\circ\CQ(T)=\id(2+q^{-4}+q^{-8})$ and
$\eps(T)=(1+q^{-2})(1+q^{-4})$ and their subtraction from $\eta'$
makes $\eta_K$ the leading term and $\theta\tens\theta$ $O(q-1)$
relative to it. This $\eta_K$ is a q-deformation of $\rho_W$ of
the usual split Casimir $X_+\tens X_-+X_-\tens X_++\h H\tens H$,
as it should be. The $\theta\tens \theta$ is a kind of `null mode'
that does not affect the geometry and could be dropped from the
metric. \eproof

\note{The computation of the moduli of torsion-free and
cotorsion-free connections for this quantum metric as well as the
Dirac operator in this case are beyond our present scope and will
be studied elsewhere. For such results we turn now to the more
tractable finite group case.}

\section{Finite Riemannian Geometry}

In this section we apply the general results above to the special
case of $M=\C[G]$ the  algebra of functions on a finite group $G$.
We first specialise the results of Section~3 to $M=\C[\Sigma]$ for
$\Sigma$ a finite set and list the main formulae of Riemannian
geometry for this case in a self-contained manner that could be put
on a computer. We then proceed to concentrate on the group case as
good source of examples where there are clear choices for the
differential structures etc. Finally, we compute everything for the
permutation group $S_3$ including solving for a canonical torsion
free cotorsion free or `Levi-Civita' spin connection in it.

\subsection{Riemannian geometry on finite sets}

Here we will see that even finite sets can be endowed with a rich
variety of `manifold' structures using the framework of Section~3.
In fact it is not true that every differential calculus on a finite
set is parallelizable (see below); i.e. there may be a still more
general theory over finite sets where we specialise the global
constructions of Section~2. This is not relevant to the finite
group case which is our main goal, and will therefore be considered
elsewhere. On the other hand, we keep the fiber of the frame bundle
to be a Hopf algebra $H$ equipped with a bicovariant differential
calculus defined by $\Omega_0$ of dimension $n$, since no special
simplification is afforded by specialising further for the tensor
product bundle. To be as concrete as possible (we have in mind
actual matrix computations for numerical gravity on finite sets)
let us assume that $H$ is finite-dimensional and choose a basis
$\{e_i\}$ for it with $e_0=1$ and $\tilde\pi_{\Omega_0}(e_i)=e_i$
for $1\le i\le n$ with the image here a basis of $\Omega_0$ (and
zero otherwise). In this way we identify $\Omega_0$ with its lift
in $H$. The dual basis $\{f^i\}$ similarly splits with $1\le i\le
n$ a basis of $\Omega_0^*$. The coproduct is of course
\eqn{Deltafin}{ \Delta e_i=c_i{}^{jk}e_j\tens e_k,} for some
structure constants.  Finally, we  write right $H$-comodules $V$
explicitly as left actions $\rho_V$ of $H^*$. We define (since we
typically convert right actions to left ones by $S^{-1}$) the
matrices
\eqn{taufin}{\tau^i=\rho_V(S^{-1}f^i)} In fact the formulae below
in the tensor product bundle depend only in this coalgebra and the
choice of quotient space (so that similar formulae hold for
coalgebra bundles\cite{BrzMa:geo} as well except that we would
specify the matrices $\tau^i$ or right action of $H^*$ directly.)

Next, we let $\Sigma$ be a finite set and $M=\C[\Sigma]$ spanned by
delta-functions  $\{\delta_x\}$ for $x\in\Sigma$. It is easy to see
(and well-known) that a general differential calculus $\Omega^1(M)$
corresponds to a subset  \eqn{1formset}{ E\subseteq
\Sigma\times\Sigma-{\rm diagonal},\quad
\Omega^1(M)=\{\delta_x\tens\delta_y|\ (x,y)\in E\}=\C E,}  where
we set to zero delta-functions corresponding to the complement of
$E$  and identify the remainder with their lifts as shown. If
$f=\sum f_x\delta_x$ is a function with components $f_x$, then
$\extd f$ has components $(\extd f)_{x,y}=f_y-f_x$ for $(x,y)\in
E$.

\begin{lemma} A $V$-bein for a finite set $\Sigma$ is a vector space
on which $H$ coacts and 1-forms
\[ E_a=\sum_{(x,y)\in  E}E_{a,x,y}\delta_x\tens\delta_y\]
for each element of a basis  $\{e_a\}_{a\in I}$ of $V$ such that
the matrices $\{E_{a,x,y}\}$ are invertible for each $x\in \Sigma$
held fixed. A necessary and sufficient condition for the existence
of a $V$-bein is that $E$ is fibred over $\Sigma$, which implies
in particular that $|E|=|\Sigma|\dim(V)$. \end{lemma} \proof We
write $E_a=e(e_a)$, etc. In principle we require the matrices
$s_e{}^{z,a}_{x,y}=\delta^z_yE_{a,x,y}$ to be invertible as maps
$M\tens V\to \Omega^1(M)$, but since they are left $M$-module maps
(or from their special form) we know that their inverses must also
be left $M$-module maps and hence of the form $s_e{}^{-1
x,y}_{z,a}=\delta^z_x E_a^{-1}{}^{x,y}$ for a collection of
matrices $E^{-1}_a{}^{x,y}$ inverse to the $E_{a,x,y}$ for each
$x$. This requires in particular that  for each $x\in\Sigma$ the
set $F_x=\{y|\ (x,y)\in E\}$ has the same size, namely the
dimension of $V$, i.e. that $E$ is a fibration over $\Sigma$ (and
$E_{a,x,y}$ is a trivialisation of the vector bundle with fiber
$\C F_x$ over $x$). The fibration is also sufficient for the
existence of a trivialisation since bundles over finite sets are
trivial. Indeed, a natural `local' class of $V$-beins is just
given by any collection of bijections $s_x:I\isom F_x$ with
$E_{a,x,y}=\delta_{s_x(a),y}$. \eproof

Similarly a $V$-cobein is a collection of 1-forms with components
$E^{*a}_{x,y}$ with respect to a dual  basis $\{f^a\}$ and with
the matrices  $\{E^{*a}_{x,y}\}$ invertible for each $y\in\Sigma$
held fixed. The metric is then \eqn{gfin}{  g=\sum_{(x,y,z)\in
F}g_{x,y,z}\delta_x\tens\delta_y\tens\delta_z,\quad
g_{x,y,z}=E^{*a}_{x,y}E_{a,y,z},\quad F=\{(x,y,z)\in \Sigma^3|\
(x,y),(y,z)\in E\},} where $\Omega^1(M)\tens_M\Omega^1(M)=\C F$.
Moreover, a connection or gauge field with values in the dual  of
$\Omega_0$ is clearly a collection of 1-forms with components
$A_{i,x,y}$. In our case $H$ coacts on $V$ so that it plays the
role of frame transformations in the frame bundle approach. In
that case $A$ induces a covariant derivative on 1-forms
\eqn{nablafin}{ (\nabla\alpha)_{x,y,z}=
(\alpha^a_y-\alpha^a_x)E_{a,y,z}-\alpha^a_xA_{a,x,y}E_{b,y,z}
\tau^i{}^b{}_a,} where $\alpha=\alpha^aE_a$ defines the  component
functions $\alpha^a$ of a 1-form $\alpha$ in the $V$-bein basis.

Next we specify $\Omega^2(M)$ by a bimodule surjection
$\wedge:\Omega^1(M)\tens_M\Omega^1(M)\to \Omega^2(M)$.

\begin{lemma} The surjections $\wedge$ are necessarily given by quotients
$V_{x,z}$ of the spaces $\C F_{x,z}$ where $F_{x,z}=\{y\in\Sigma|\
(x,y,z)\in F\}$ such that the image of the vector $(1,1,\cdots,
1)$ is zero whenever $(x,z)\notin E$ with $x\ne z$. Explicitly,
\eqn{wedgefin}{(\wedge f)_{x,\alpha,z}=\sum_{y\in
F_{x,z}}f_{x,y,z}p_{x,z}{}^y{}_\alpha,\quad } for a family of
matrices $p_{x,z}$ with respect to a basis $\{e_\alpha\}$ of each
$V_{x,z}$ and with rows summing to zero when $(x,z)\notin E$ with
$x\ne z$.
\end{lemma}
\proof We require to quotient $\Omega^1(M)\tens_M\Omega^1(M)=\C F$
by a subbimodule. This must therefore take the form shown for some
surjections $p_{x,z}$. The additional stated condition is for
$\extd^2=0$ (so the maximal prolongation will be $V_{x,z}=\C
F_{x,z}$ when $(x,y)\in E$ and $\C F_{x,z}/\C.(1,1,\cdots,1)$
otherwise). The argument is similar to that in \cite{BrzMa:dif}.
There may be additional restrictions imposed by requiring the
$\Omega^2(M)$ to be part of a global $\Omega^2(P)$ as explained in
Section~3. \eproof

When $\alpha_{x,y},\beta_{x,y}$ are the components of 1-forms as
above then \eqn{extd2fin}{ (\extd\alpha)_{x,\alpha,z}=\sum_{y\in
F_{x,z}}
(\alpha_{x,y}+\alpha_{y,z}-\alpha_{x,z})p_{x,z}{}^y{}_\alpha,\quad
(\alpha\wedge\beta)_{x,\alpha,z}=\sum_{y\in
F_{x,z}}\alpha_{x,y}\beta_{y,z}p_{x,z}{}^y{}_\alpha.} With such an
explicit description of $\Omega^2(M)$ it is clear that a
connection $A$ is regular if \eqn{regfin}{ \sum_{1\le j,k\le
n,y}c_i{}^{jk} A_{j,x,y}A_{k,y,z}p_{x,z}{}^y{}_\alpha=0,\quad
\forall q\notin\CC\cup\{e\}.}   Its  curvature is  \eqn{Ffin}{
F_{i,x,\alpha,z}=(\extd A_i)_{x,\alpha,z}+\sum_{1\le j,k\le
n,y}c_i{}^{jk}A_{j,x,y}A_{k,y,z}p_{x,z}{}^y{}_\alpha.} The actual
Riemann tensor is the 2-form valued operator on 1-forms,
\eqn{riemfin}{
R_{x,\alpha,z}{}^a{}_b=F_{i,x,\alpha,z}\tau^i{}^a{}_b,\quad
R\alpha=\alpha^a R^b{}_a\tens_M E_b.} Meanwhile, the zero torsion
and  zero cotorsion equations are vanishing of \eqn{torfin}{(\bar
D \wedge e)_{a,x,\alpha,z}=(\extd
E_a)_{x,\alpha,z}+\sum_{i,b,y}A_{i,x,y}E_{b,y,z}
p_{x,z}{}^y{}_\alpha\tau^i{}^b{}_a,} \eqn{cotorfin}{(D\wedge
e^*)^a_{x,\alpha,z}=(\extd
E^{*a})_{x,\alpha,z}+\sum_{i,b,y}E^{*b}_{x,y}A_{i,y,z}
p_{x,z}^y{}_\alpha\tau^i{}^a{}_b.}

Also, a `lift' $i:\Omega^2(M)\to \Omega^1(M)\tens_M\Omega^1(M)$ is
given similarly to the discussion above by a collection of
inclusions $i_{x,z}:V_{x,z}\to \C F_{x,z}$ or a family of
rectangular matrices $i_{x,z}{}^\alpha{}_y$. We let
$\pi_{x,z}=i_{x,z}\circ p_{x,z}$ so that $\pi^y{}_w=p^y{}_\alpha
i^\alpha{}_w$ at each $x,z$. If $i$ is a true lift so that $p\circ
i=\id$ then $i\circ\wedge$ is a projection splitting
$\Omega^1(M)\tens_M\Omega^1(M)$ into something isomorphic to
$\Omega^2(M)$ plus a complement and the $\pi_{x,z}$ are likewise a
family of projection matrices. We do not want to strictly assume
this, however. Given $i$, we have an interior product and, in
particular, a Ricci tensor  \eqn{riccifin}{ {\rm
Ricci}_{x,y,z}=i(F_i)_x^{ab}E_{b,x,y}E_{c,y,z}\tau^i{}^c{}_a.}
Here $i(F_i)_{x,w,z}$ in $\Omega^1(M)\tens_M\Omega^1(M)$ is as in
(\ref{Ffin}) but with $\pi_{x,z}{}^y{}_\alpha$ in place of
$p_{x,z}{}^y{}_\alpha$ written there. The $V$-bein components here
are defined by $i(F_i)_{x,y,z}=i(F_i)_x^{ab}E_{a,x,y}E_{b,y,z}$ as
usual.

Finally, gamma-matrices are a collection of matrices  $\gamma_a$
acting on spinors $\psi$ which are functions with values in a
vector space $W$ on which $H$ coacts by $\rho_W$, say. We define
the corresponding matrices $\tau^a_W$ as above. Then the
associated Dirac  operator is \eqn{dirfin}{
\Dsl=\del^a\gamma_a-A_i^a\gamma_a \tau^i_W.}

For the case when $H=\C[G]$ it is actually useful to chose a
different basis for $H$ that reflects better the group structure,
namely we label the basis by the group elements themselves (so
$e_i$ is the delta-function at $i\in G$ and
$c_i{}^{jk}=\delta_i^{jk}$). This has the same form as above
except that the old $e_0$ above is the sum of all the new basis
elements. The role of $e_1,\cdots,e_n$ is played by $e_i$ for
$i\in\CC$ a subset of order $n$ not containing the identity
element $e\in G$ (see below), which is purely a notational change.
All the formulae above have the same form in this case except the
regularity and curvature equations, for which one has to make a
careful change of basis (or use the form of $\tilde\pi_{\Omega_0}$
in the new basis as given in the next section). One has instead,
\eqn{regfinG}{
\sum_{jk=q,y}A_{j,x,y}A_{k,y,z}p_{x,z}^y{}_\alpha=0,\quad \forall
q\notin\CC\cup\{e\}}\eqn{FfinG}{ F_{i,x,\alpha,z}=(\extd
A_i)_{x,\alpha,z}+\sum_{jk=i,y}A_{j,x,y}A_{k,y,z}p_{x,z}^y{}_\alpha-
\sum_{j,y}(A_{j,x,y}A_{i,y,z}+A_{i,x,y}A_{j,y,z})p_{x,z}^y{}_\alpha}
respectively, with $i,j,k\in\CC$.

Whereas the above tensorial formulae are suitable for numerical
computations, let us note finally that we also have more algebraic
`Cartan calculus' formulae based on (\ref{partialM}).  Thus,
\eqn{setcartan}{ E_a f=\rho_a{}^b(f)E_b,\quad \extd
f=[\theta,f];\quad \rho_a{}^b(f)(x)=\sum_{y\in F_x}
E_b^{-1xy}f(y)E_{axy},\quad \theta=\theta^aE_a,\quad
\theta_a(x)=\sum_{y\in F_x}E_a^{-1xy}} for all functions $f$. For
$\Omega^2(M)$ we can build $\wedge$ from a $G$-equivariant
projector $\pi(e_a\tens e_b)\equiv \pi_{ab}{}^{cd}e_c\tens e_d$ on
$V\tens V$. Then \eqn{setpi}{
\pi_{x,z}{}^y{}_w=\pi_{ab}{}^{cd}E_a^{-1xy}E_b^{-1yz}E_{cxw}E_{dwz}}
for the above family of projection matrices. This imposes
contraints on $(\pi,E)$ and defines a moduli space of
$G$-parallelizable manifold structures on a finite set of a given
order.

\subsection{Riemannian geometry on finite groups}

We now specialise further to the case the case $M=H=\C[G]$ with
the same bicovariant differential calculus on both. This gives a
nontrivial setting at the level of finite groups. In principle one
obtains `geometric invariants' of finite groups equipped with a
differential calculus (i.e. a conjugacy class), which is certainly
of independent mathematical interest as well as of physical
interest as a simple toy setting for finite gravity.

As mentioned above, the coirreducible calculi are classified
immediately from \cite{Wor:dif} by   nontrivial conjugacy classes
$\CC\subset G$. In fact we do not need  to assume that the
calculus is irreducible and hence in what follows $\CC$ is any
$\Ad$-stable subset not containing the group unit element $e\in
G$. We denote the elements of $\CC$ by $a,b,c$, etc. Then
\eqn{QG}{ Q_H=\{\delta_q|\ q\ne e,\ q\notin \CC\},\quad
\Omega_0=\{\delta_a|\ a\in \CC\}=\C\CC,\quad
\Ad(\delta_a)=\sum_{g\in G}\delta_{gag^{-1}}\tens \delta_g}
\eqn{extdG}{ \extd f=\sum_{a\in\CC} (\del^a f)\cdot
E_a,\quad\del^a=R_a-\id,\quad E_a\cdot f=R_a(f)\cdot E_a} where
$R_a(f)=f((\ )a)$. In this description we identify a basis element
$e_a$ of $\Omega_0$ with a fixed lift $\delta_a\in\ker\eps$, which
is an $\Ad$-invariant identification. The projection from $\C[G]$
to $\Omega_0$ is then \eqn{projG}{\tilde\pi_{\Omega_0}(\delta_g)=
\begin{cases}\delta_g&\text{if
$g\in\CC$}\\ -\sum_{a\in\CC}\delta_a&\text{if $g=e$}\\
0&\text{else.}\end{cases}} The elements of $\Omega_0$ viewed in
$\Omega^1(H)$ are the values of the Maurer-Cartan form
$e:\Omega_0\to \Omega^1(H)$, \eqn{mcformG}{
E_a=e(\delta_a)=\pi_{N_H}(\sum_{g\in G}\delta_{g}\tens
\delta_{ga}),\quad N_H=\{\delta_g\tens\delta_{gq}|\ g\in G,\ q\ne
e,\ q\notin\CC\}.} In terms of the general finite set case, we
have a local form of the $V$-bein and the action, \eqn{finG}{
s_x(a)=xa,\quad E_{a,x,y}=\delta_{xa,y},\quad
\tau^a{}^b{}_c=\delta^b_{a^{-1}ca}-\delta^b_c.}

The rest of our treatment in the finite group case, is more easily
handled in the `Cartan calculus' form at the end of Section~5.1
i.e. by algebraic relations among the $\{E_a\}$ generators of the
entire exterior algebra rather than in `spacetime coordinates'
$\alpha_{x,y}$, etc. Thus, the higher exterior algebra is
generated\cite{Wor:dif} by the relations at the first order and
the additional relations implied by the braiding \eqn{braidG}{
\Psi(E_a\tens_H E_b)=E_{aba^{-1}}\tens_H E_a.} Thus, in
$\Omega^2(H)$ the quotient to the wedge product consists in
setting to zero all  linear combinations invariant under $\Psi$.
In particular, for all $g\in G$ the elements $\sum_{ab=g} E_a\tens
E_b$ are invariant (after a change of variables), hence these
along with the clearly invariant $E_a\tens E_a$ give some
immediate relations \eqn{rel2G}{ \sum_{a,b\in \CC,\ ab=g}
E_a\wedge E_b=0,\quad \forall g\in G,\quad E_a\wedge E_a=0.} Using
these relations and (\ref{projG}) the Maurer-Cartan equation on
any Hopf algebra becomes \eqn{mcG}{ \extd E_a-\sum_b (E_a\wedge
E_b+E_b\wedge E_a)=0.} Meanwhile, the partial derivatives
trivially obey \eqn{delG}{ \del^a(mn)=m\del^an + (\del^a
m)R_a(n),\quad \del^a\del^b=\del^{ab}-\del^a-\del^b} as
$R_{ab}=R_a R_b$, where we extend the same definitions to $ab\in
G$. One can also write
\eqn{worlieG}{\del^a\del^b-\del^b\del^{b^{-1}ab}=\del^{b^{-1}ab}-\del^a}
as some form of Lie algebra\cite{Wor:dif}, however such a point of
view can only be taken so far, and we do not use it. Rather, the
$\del^a$ form a representation of a braided-Lie
algebra\cite{Ma:lie}. The above formulae, with the exception of
our notations such as (\ref{projG}) and the observation
(\ref{rel2G}), are all immediate from the general theory of
\cite{Wor:dif} and are the starting point of any quantum-groups
inspired  noncommutative geometry on finite groups. We will also
need a more full description of $\Omega^2(H)$ in the finite group
case, provided by the following lemma.

\begin{lemma} For all $g\in G$ let $P_g=(\C \CC\cap g\CC^{-1})^{\sigma}$ the
invariant subspace of the vector space with basis $\CC\cap
g\CC^{-1}$, where $\sigma$ sends a basis element $a$ to $a^{-1}g$.
Let $\{\lambda^{g,\alpha}\}$ be a basis of $P_g$. Then the
relations of $\Omega^2(H)$ are
\[ \sum_{a,b\in \CC,\ ab=g} \lambda^{g,\alpha}_a E_a\wedge E_b=0,\quad \forall
g\in G,\quad \forall\alpha.\]
\end{lemma}
\proof For any $\lambda\in P_g$, we clearly have invariance under
$\Psi$ as $\sum_{ab=g}\lambda_a E_{aba^{-1}}\tens_H
E_a=\sum_{cd=g}\lambda_{c^{-1}g}E_c\tens_H E_d=\sum_{ab=g}\lambda_a
E_a\tens_H E_b$, by the $\sigma$-invariance of $\lambda$. Hence
relations of the form shown hold in $\Omega^2(H)$ for any basis of
$P_g$, for each $g\in G$. One can show that this is a full set of
relations after a detailed analysis of the kernel of $\id-\Psi$ in
this case.
\eproof

We are now ready to specialise our results of Sections~3,4 to
obtain a theory of Riemannian geometry for finite groups. First of
all, as a trivial example of the theory in \cite{Ma:lie} we may
view $\Omega^*_0$ as the image under $\pi_\eps$ of a braided-Lie
algebra with trivial background braiding and \eqn{LG}{L=\{x^a|\
a\in\CC\}=\C\CC, \quad [x^a,x^b]=x^{b^{-1}ab},\quad \Delta
x^a=x^a\tens x^a.} The braided enveloping bialgebra $U(L)$ from
\cite{Ma:lie} in this case (because the background braiding is
trivial) is actually a usual bialgebra or quantum group without
antipode. It comes with a bialgebra homomorphism to the group
algebra $\C G$, \eqn{U(L)G}{ U(L)=\C\<x^a\>/\ x^a
x^b=x^bx^{b^{-1}ab},\quad p:U(L)\to \C G,\quad p(x^a)=a.} The
further projection of the braided-Lie algebra generators to the
kernel of the counit gives the basis $\{f^a=a-e\}$ of $\Omega_0^*$
dual to the $e_a=\delta_a$ via Hopf algebra duality. One may also
consider `braided gauge theory' with $A$ having values in $L$
rather than in $\Omega_0^*$, but for the present we need this
theory mainly to have a braided-Killing form.

\begin{propos}  The braided-Killing form of $L$ is a symmetric
positive-integer valued and $\Ad$-invariant bilinear form on the
conjugacy class given by \[ \eta(x^a,x^b)=n(ab)\equiv \#\{c\in
\CC|\ c ab=ab c\}=\eta(x^b,x^{b^{-1}ab}),\quad \forall
a,b\in\CC.\] We say that a conjugacy class is {\em semisimple} if
this associated Killing form is nondegenerate.
\end{propos} \proof The braided-Killing form is defined as the trace
$\eta(x^a,x^b)=\sum_{c\in\CC}\<\delta_c,[[x^c,x^a],x^b]\>$ which is
clearly as shown (the number of $c\in\CC$ commuting with $ab$). Its
formal properties are part of the general theory of braided-Lie
algebras. The relevant braiding in the present case is that of the
category of crossed $\C G$-modules and has the form $\Psi(x^a\tens
x^b)=x^b\tens x^{b^{-1}ab}$ (as above for differentials), so that
$\eta$, depending only on the product, is clearly braided-symmetric
in the sense $\eta=\eta\circ\Psi$. Hence it is also symmetric in
the usual sense (because $S^2=\id$ in Proposition~4.2.)
Ad-invariance is clear as well. Note that the braided-Lie algebra
itself is bosonic as the category of $\C[G]$-comodules in which it
lives has a trivial background braiding. \eproof

Because $\Omega_0^*$ can be identified naturally (and
$\Ad$-invariantly) with $L$ viewed inside $\C G$, we pull back
this braided-Killing form to obtain on $\Ad$-invariant bilinear
form \eqn{etafxG}{\eta(f^a,f^b)\equiv \eta(x^a,x^b)=n(ab).} Note
that this gives a slightly different Killing form than the trace
of $\Ad_{f^a}\Ad_{f^b}$ i.e. taking `Lie bracket' as the quantum
group adjoint action of $f^a=a-e$ in $\C G$, giving instead
\[\eta(f^a,f^b)=n(ab)-n(a)-n(b)+n(e)\] more similar to Theorem~4.7.
This is also $\Ad$-invariant so (if nondegenerate) could also be
used to define a metric with essentially the same geometry.  Also
note that for an Abelian group the conjugacy classes are
singletons but we make take  $\CC$ a collection of these and the
same formulae as above for a (reducible) differential calculus.
The Killing form will be degenerate in this case but the
$\delta$-function provides instead a suitable symmetric and
invariant bilinear form.  Let \eqn{etaabG}{
\eta^{ab}=\eta(f^a,f^b),\quad
\eta^{-1}_{ab}=\eta^{-1}(\delta_a,\delta_b)} for the braided
Killing form and its inverse in our basis. We will write
$\alpha=\alpha^a\cdot E_a$ for any 1-form and similarly for the
components of  higher cotensors. We sum over repeated indices
$a,b\in\CC$ in tensor expressions.

In this setting the main equations of `quantum group Riemannian
geometry' of Section~3 become as follows. A framing is a
collection of 1-forms $\{E_a\}$ such that every 1-form is a unique
linear combination of these with coefficients functions from the
left (e.g. as above). A spin connection is a collection $\{A_a\}$
of 1-forms and the covariant derivative associated to a spin
connection and framing on any 1-form $\alpha$ is
\eqn{covG}{\nabla\alpha=\extd\alpha^a\tens_H E_a-\alpha^a\sum_b
A_b\tens_H (E_{b^{-1}ab}-E_a).} The extra term in $\nabla$ comes
from the projection $\tilde\pi_{\Omega_0}$ or equivalently from
the fact that the role of `Lie algebra' is being played by the
vectors $\{a-e\}$ in the group algebra as explained above. The
associated torsion tensor $T:\Omega^1(H)\to\Omega^2(H)$ measures
the deviation $T\alpha=\extd\wedge\alpha-\nabla\alpha$ and the
zero-torsion condition is vanishing of \eqn{torG}{ \bar D_A\wedge
E_a=\extd E_a+\sum_b A_b\wedge (E_{b^{-1}ab}-E_a).} The curvature
$\nabla^2$ associated to a regular connection corresponds to a
collection of two forms $\{F_a\}$ defined by \eqn{curvG}{
F_a=\extd A_a+\sum_{cd=a,\ c,d\in\CC}A_c\wedge A_d-\sum_b
(A_b\wedge A_a+A_a\wedge A_b)} where regularity in Section~3
becomes the condition \eqn{regG}{ \sum_{ab=q,\ a,b\in\CC}
A_a\wedge A_b=0,\quad \forall q\ne e,\ q\notin\CC.} It ensures
that the curvature descends to $\Omega_0$, otherwise it
potentially has values $\{F_g\}$ for all $g\ne e$. It is clear
that the Maurer-Cartan form can be viewed as a regular connection
with zero curvature. For any connection the associated Riemann
curvature is the 2-form-valued operator
\eqn{rieG}{R\alpha=\alpha^a\sum_b F_b\tens_H (E_{b^{-1}ab}-E_a)}
on 1-form $\alpha=\alpha^aE_a$, according to the correspondence in
Section~3.

To define the Ricci tensor (or to define interior products in
general) we need a bimodule inclusion or `lift' $\Omega^2(H)\to
\Omega^1(H)\tens_H\Omega^1(H)$. The obvious one for the
bicovariant calculus, although not precisely a lift any more (not
covered by $\wedge$) is provided by  \eqn{worliftG}{
i=\id-\Psi,\quad i(E_a\wedge E_b)=E_a\tens_H
E_b-E_{aba^{-1}}\tens_H E_a.} We provide now another possibility
which {\em is} actually a lift in a natural manner, so that
$i\circ\wedge$ is an actual projection operator on
$\Omega^1(H)\tens_H\Omega^1(H)$.

\begin{propos} For $H=\C[G]$ and the bicovariant calculus, there is a
canonical splitting of $\wedge$ to a bimodule projection operator,
defined by
\[ i(E_a\tens_H E_b)=E_a\tens_H
E_b-\sum_{\alpha}\mu^{\alpha,a}\sum_{cd=g}\lambda^{\alpha}_cE_c\tens_H
E_d\] where $\mu^{\alpha}\in P_g$ are a dual basis to the
$\lambda^{\alpha}$ with respect to the dot product as vectors in
$\C \CC\cap g\CC^{-1}$, and $g=ab$ is fixed. \end{propos} \proof
Here the summed terms vanish under $\wedge$ by Proposition~5.1, so
that $i$ as stated indeed splits this for any choice of
coefficients $\mu^{\alpha}$.  We choose these to be the dual basis
to the $\lambda^{\alpha}$ so that $\sum_{a\in\CC\cap
g\CC^{-1}}\mu^{\alpha,a}\lambda^\beta_a=\delta^{\alpha,\beta}$.
Here $g=ab$ is suppressed in our notation. Then one may verify that
the map is well-defined on $\Omega^2(H)$, i.e.
$\sum_{ab=g}\lambda_a i(E_a\wedge E_b)=0$. Finally, $i$ by
definition extends as a left $H$-module map and, since we only add
terms of the same `total degree' $g$ with respect to the right
action, it becomes also a right module map. \eproof

Given the choice of `lift', the Ricci tensor constructed from the
Riemann tensor by making a point-wise trace over the input and the
first output of $i(R)$, is \eqn{ricciG}{{\rm Ricci}=
\sum_{a,b,c}i(F)^{ab}_c E_b\tens_H (E_{c^{-1}ac}-E_a),} where
$i(F_c)=i(F)^{ab}_c E_a\tens E_b$.

Next, a gamma-matrix is a collection of endomorphisms
$\{\gamma_a\}$ of a vector space $W$ on which $G$ acts by a
representation $\rho_W$ say, subject to further constraints to be
discussed on the $\gamma$. A `spinor' field is a $W$-valued
function on $G$, and \eqn{diracG}{
\Dsl=\del^a\gamma_a-A_b{}^a\gamma_a\tau_W^b,\quad
\tau^a=\rho_W(a^{-1}-e)} where $A_b=A_b{}^aE_a$ determines the
components of each $A_a$. The $\tau_W^a$ are the `Lie algebra'
generators $f^a$ in the representation $\rho_W$. The group inverse
here makes them actually a right-action rather than a left one
(just as the $\del^a$ are actually right-derivations).

Finally, a metric is determined by a choice of framing and a
coframing $\{E^*{}^a\}$ which is a collection of 1-forms such that
every 1-form is a unique combination of these with coefficient
functions from the right. Given a framing, a general coframimg and
hence a general metric is determined by a point-wise invertible
function-valued matrix $\{g_{ab}\}$ and given as a cotensor by
\eqn{metricG}{  g=E_a g^{ab}\tens_H E_b} where $g^{ab}$ is the
matrix inverse (e.g. $g^{ab}=\eta^{ab}$ above).  The cotorsion of
the  spin connection is the torsion with respect to the coframing
and corresponds to \eqn{cotorG}{ D_A E^*{}^a=\extd
E^*{}^a+\sum_c(E^{*cac^{-1}}-E^{*a}) \wedge A_c.} Vanishing of the
cotorsion generalises the notion of metric compatibility in a
slightly weaker `skew' formulation appropriate to our not
requiring the metric symmetric\cite{Ma:rie}.    One is at liberty
now to do `finite gravity'. That is one can look at the moduli
spaces for the above data and solve the various equations as well
as others such as given by the variation or minimisation of an
action. The role of Einstein-Hilbert action can be played for
example by the trace of $\Dsl^2$. Since everything is finite we do
not need to worry about regularisations and Dixmier traces etc. as
in the approach of Connes\cite{Con:geo}. We will not attempt this
here but we will show for a nontrivial example in the next section
that the moduli space of our basic data is not empty.

For example, we could  fix the framing and coframing to be the
natural ones on any quantum group defined as above by the
Maurer-Cartan form and a `braided-Killing form' $g=\eta$. We have
established these canonical choices in Section~4. If one wants a
torsion free spin connection we then we have to solve (in view of
the Maurer-Cartan equations already obeyed), the condition
\eqn{torcanG}{ \sum_{b\ne a} A_b\wedge(E_{b^{-1}ab}-E_a)+E_b\wedge
E_a+E_a\wedge E_b=0, \quad \forall a\in\CC.} We need only solve
this for all except one $a$ since the sum over $a$ is
automatically zero in view of (\ref{rel2G}). Finally we could take
for $\gamma_a$ the `tautological' one in Section~4,
\eqn{cangammaG}{\gamma_a=\sum_b\eta^{-1}_{ab}\rho_W(b-e).} These
are equivariant and obey \eqn{cangammanorm}{
\eta^{ab}\gamma_a\gamma_b=\rho_W(C),\quad
C=\sum_{a,b\in\CC}\eta^{-1}_{ab}(a-e)(b-e)} where $C$ is the
braided Casimir element associated to the braided-Killing form.

We can also consider other choices of gamma-matrices
$\{\gamma_a\}$. Our other new proposal mentioned in Section~4 is
that the gamma-matrices could be restricted by the requirement
that Connes prescription\cite{Con:geo} for the exterior algebra
$\Omega^\cdot_{\Dsl}$  obtained from $\Dsl$ should coincide with
our bicovariant approach above, which would be the case
classically. This condition is independent of the choice of
framing, coframing or spin-connection since the commutators
$[\Dsl,m]$ relevant for this ($m$ any function) are independent of
these. The following proposition shows, however, that this is {\em
not} necessarily a natural restriction  in the present context of
finite groups.

\begin{theorem} A necessary condition for the Connes exterior
algebra induced by $\Dsl$ to contain the relations of the
Woronowicz bicovariant one on $\C[G]$ is \[ \gamma_a^2=0,\quad
{\rm if}\ a^2\in\CC\cup\{e\},\ a\in\CC,\quad \sum_{ab=g,\
a,b\in\CC}\gamma_a\gamma_b=0,\quad \forall g\in\CC\cup\{e\}.\]
\end{theorem} \proof We recall that \cite{Con:geo} considers a
representation $\pi_\Dsl$ of  the universal exterior algebra a
spectral triple. The relevant part of this construction, however,
does not depend on Hilbert spaces or self-adjointness and works
for any algebra $M$ and operator $\Dsl$ on a vector space in which
$M$ is also represented. In our case the algebra is $M=H=\C[G]$
and the vector space is of the form $M\tens W$ and $M$ is
represented by multiplication. Then \[ \pi_\Dsl:\Omega^\cdot M\to
\End(M\tens W),\quad \pi_\Dsl(m\tens n\tens \cdots\tens
p)=m[\Dsl,n]\cdots [\Dsl,p]\] defines the exterior algebra
$\Omega^\cdot_\Dsl$ as the quotient of the universal one modulo
the differential graded ideal generated by the kernel of
$\pi_\Dsl$. At degree 1 we  know from Section~3 that \[
m[\Dsl,n]=\sum_{a\in \CC}m(\del^a n)R_a\tens\gamma_a\] from which
it is clear that for an injective map $\gamma:\Omega_0\to\End(W)$
the kernel of $\pi_\Dsl$ at degree 1 is the same as $N_H$, the
ideal set to zero by $m\extd n=m(\del^a n)E_a$. At degree 2 we
have \[ [\Dsl,m][\Dsl,n]=\sum_{c,d\in \CC}(\del^c m)\circ
R_c\circ(\del^d n) R_d\tens \gamma_c\gamma_d =\sum_{c,d\in
\CC}(\del^c m)(\del^d R_c n)R_{dc}\tens\gamma_c\gamma_{c^{-1}dc}\]
after a change of variables. Next, working in the universal
calculus, the product of Maurer-Cartan forms is \[
e(\delta_g)\tens_H e(\delta_h)=\sum_{b\in
G}\delta_{b}\tens\delta_{bg}\tens\delta_{bgh}\] and one finds
\begin{eqnarray*}\pi_\Dsl( e(\delta_g)\tens_H e(\delta_h))&&=\sum_{c,d\in\CC,
b\in
G}\delta_{b}(\delta_{bgc^{-1}}-\delta_{bg})(\delta_{bghc^{-1}d^{-1}}
-\delta_{bghc^{-1}})R_{dc}\tens \gamma_c\gamma_{c^{-1}dc}\\
&&=\sum_b\delta_{b}R_{gh}\tens\gamma_g\gamma_h=R_{gh}\tens\gamma_g\gamma_h
\end{eqnarray*} where only the leading term in each difference contributes
when $g\ne e$ and $h\ne e$. The $\delta$-functions then fix $c=g$
and $d=ghg^{-1}$ provided $h\in\CC$ (otherwise we obtain zero).
Also,  \begin{eqnarray*} \pi_\Dsl(\extd N_H)&&=\pi_\Dsl\{(\extd
\delta_g)\tens_H e(\delta_q)+\delta_g\extd
e(q)|\ g\in G,\ q\ne e,\ q\notin\CC\}\\
&&=\{\delta_gR_q\tens\sum_{c,d\in\CC,
dc=q}\gamma_c\gamma_{c^{-1}dc}|\ g\in G,\ q\ne e,\ q\notin\CC\}\\
&&=\{\delta_gR_q\tens\sum_{a,b\in\CC, ab=q}\gamma_a\gamma_b|\ g\in
G,\ q\ne e,\ q\notin\CC\}\end{eqnarray*} where the first term
fails to contribute since it is of the form a function times
$e(\delta_h)\tens_H e(\delta_q)$ where $h,q\ne e$ and $q\notin\CC$
(see above). The second term only contributes
$\pi_\Dsl(\delta_g\tens\delta_{bh^{-1}}\tens\delta_b)$ which comes
out as stated by similar computations to those above and a further
change of variables as shown. Hence $\pi_\Dsl$ applied to the
expressions leading to the relations   (\ref{rel2G}) in the
Woronowicz calculus, namely \[ R_{g}\tens
\sum_{ab=g}\gamma_a\gamma_b,\quad R_{a^2}\tens \gamma_a^2\] do not
lie in $\pi_\Dsl(\extd N_H)$ if $g,a^2\in\CC\cup\{e\}$
respectively, unless zero. Hence for the Woronowicz ideal at
degree 2 to be contained in $\ker\pi_\Dsl+\extd N_H$ a necessary
condition is for these operators to vanish when
$g,a^2\in\CC\cup\{e\}$. This gives the conditions stated. \eproof

This will often be a sufficient condition as well, for suitably
non-degenerate $\gamma$ and in the nice cases where (\ref{rel2G})
are all the relations (at least at degree 2). Moreover,  when the
conclusion holds it often means that the Connes and Woronowicz
calculi actually coincide, because the Woronowicz one tends to
have the most relations anyway in practice.  The theorem is a
surprising result but easily verified for example on $\C[\Z_2]$.
This has only one nontrivial conjugacy class $\CC=\{u\}$ where $u$
with $u^2=e$ is the nontrivial element of $\Z_2$. The Woronowicz
calculus has $E_u\wedge E_u=0$ and hence $\Omega^2=0$ while the
Connes prescription can give this (if) and only if $\gamma_u^2=0$.
The nilpotency is associated to the order 2 of $u$ and means in
particular that the Dirac operator itself will not typically be
Hermitian with respect to the obvious inner products. Such
nilpotent models could still be physically interesting and one of
them, on $\C[\Z_2\times\Z_2]$, will be explored elsewhere as a
model where Connes' approach and the quantum groups approach to
the discrete part of the geometry intersect\cite{MaSch:lat}. One
may easily make the same analysis in the general setting of
Section~4 for any Hopf algebra but this simple example is enough
to show the limitations of this approach (therefore we have
omitted the full analysis). The result means that for $\gamma$
chosen according to other criteria (such as equivariance) one will
typically have a different induced higher order calculus
$\Omega_{\Dsl}$ than the usual bicovariant one of Woronowicz
natural in this context. One may work with either one or with the
maximal prolongation with the difference appearing at $\Omega^2$
and higher, i.e. affecting the curvature and vanishing of torsion
etc.

\subsection{Riemannian geometry of $S_3$}

We now turn to a concrete example, the permutation group $G=S_3$
generated by $u,v$ with relations \eqn{S3}{ u^2=v^2=e,\quad
uvu=vuv.} The conjugacy class $\CC=\{u,v,uvu\}$ is semisimple in
the sense of Proposition~5.1 while the other nontrivial conjugacy
class $\{uv,vu\}$ is not. We therefore fix this first case, i.e.
work with a 3-dimensional bicovariant differential calculus. In
this case one finds by enumeration that
\eqn{etaS3}{\eta^{ab}=3\delta^{ab}.} The braided-Lie algebra here
is \eqn{LS3}{ [x^u,x^v]=x^w=[x^v,x^u],\quad
[x^u,x^w]=x^v=[x^w,x^u],\quad [x^v,x^w]=x^u=[x^w,x^v]} and $U(L)$
is generated by $1$ and $x^a$ with the relations \eqn{U(L)S3}{
x^ux^v=x^vx^w=x^wx^u,\quad x^vx^u=x^ux^w=x^wx^v.} If one defined
the Killing form by the adjoint action of the $f^a$ then one would
have instead $\eta{}^{ab}=3\delta^{ab}+3$. In fact any constant
offset here not change anything in terms of the resulting
connection etc. (but could render $\eta$ degenerate). The various
metrics just differ by a multiple of $\sum_{a,b}E_a\tens_H E_b$
which will turn out to play a somewhat neutral role.

The explicit form of the higher differential calculus is
well-known and in this case (\ref{rel2G}) give all the relations
at degree 2, namely \[ E_u\wedge E_u=E_v\wedge E_v=E_{uvu}\wedge
E_{uvu}=0\] \eqn{rel2S3}{E_u\wedge E_v+E_v\wedge
E_{uvu}+E_{uvu}\wedge E_u=0,\quad E_v\wedge E_u+E_{uvu}\wedge
E_v+E_u\wedge E_{uvu}=0} so that $\Omega^2(\C[S_3])$ is
4-dimensional. Lemma~5.1 establishes that these are in fact a full
set of relations in this case. The Maurer-Cartan equations
(\ref{mcG})  immediately become
\[ \extd E_u+E_{uvu}\wedge E_v+E_v\wedge
E_{uvu}=0,\quad \extd E_v+E_u\wedge E_{uvu}+E_{uvu}\wedge E_u=0\]
\eqn{mcS3}{ \extd E_{uvu}+E_v\wedge E_u+E_u\wedge E_v=0.} This has
been observed by many authors using Woronowicz bicovariant
calculus. With this background we now construct explicit solutions
to our torsion and cotorsion conditions.

\begin{propos} For the framing by the Maurer-Cartan
form, the moduli space of zero-torsion spin connections is
12-dimensional and takes the form
\[ A_u=(\alpha+1) E_u+\gamma E_v +\beta E_{uvu},\quad
A_v=\gamma E_u +(\beta+1)E_v+\alpha E_{uvu}\] \[  A_{uvu}=\beta
E_u+\alpha E_v+(\gamma+1)E_{uvu},\quad\quad
\alpha+\beta+\gamma=-1,\] where $\alpha,\beta,\gamma$ are
functions subject to the constraint shown. They obey $\sum_a
A_a=0$.
\end{propos}
\proof We solve the two equations \[ A_v\wedge
(E_{uvu}-E_u)+A_{uvu}\wedge (E_v-E_u)=E_{uvu}\wedge E_v+E_v\wedge
E_{uvu}\] \[ A_u\wedge(E_{uvu}-E_v)+A_{uvu}\wedge
(E_u-E_v)=E_u\wedge E_{uvu}+E_{uvu}\wedge E_u\] where the third
equation in (\ref{torcanG}) will be automatic given the other two.
The right hand sides here are $-\extd E_u$ and $-\extd E_v$
respectively. Into these equations we write the component
decomposition $A_a=A_a{}^b E_b$ with \[ A_u{}^u=\alpha+1,\quad
A_v{}^v=\beta+1,\quad A_{uvu}{}^{uvu}=\gamma+1\] say (these could
be functions on the group, not numbers). We then write everything
in terms of any four linearly independent 2-forms, say $E_u\wedge
E_v, E_v\wedge E_u, E_u\wedge E_{uvu}$ and $E_{uvu}\wedge E_u$,
writing the other two in terms of these via the above relations in
$\Omega^2$. The coefficients of these four 2-forms must separately
vanish and give us the four equations \[ A_{uvu}{}^u-\beta=0,\quad
-A_{uvu}{}^v-\gamma-(\beta+1)=0,\quad A_v{}^u-\gamma=0,\quad
-A_v{}^{uvu}-(\gamma+1)-\beta=0\] respectively. Similarly for the
other equation to give the solution stated.  \eproof

Next we consider metrics. The general moduli space of all metrics
is clearly $GL_3$ raised to the 6th power, as we have a reference
metric $\eta^{ab}$ provided by the braided-Killing form, and hence
a natural reference coframing $E^{*a}=\eta^{ba}E_b$. The
corresponding metric induced by the braided Killing form is of
course \eqn{getaS3}{g=\eta^{ab}E_a\tens_H E_b=3\sum_{a}E_a\tens_H
E_a.}

\begin{corol} (i) The moduli space of cotorsion-free connections
with respect to  the coframing defined by the braided-Killing form
metric $\eta^{ab}$ is also 12-dimensional and has a similar form
to the above, with coefficients $\alpha,\beta,\gamma$ on the
right. (ii) The moduli of torsion free and cotorsion free
connections is 2-dimensional, with $\alpha,\beta,\gamma$ numbers.
(iii) The point $\alpha=\beta=\gamma=-{1\over 3}$ in this moduli
space is the unique regular torsion-free and cotorsion-free or
`Levi-Civita' connection on $S_3$. This and its nonzero curvature
are \[ A_a=E_a-{1\over 3}\theta,\quad \theta=\sum_a E_a,\quad
F_a=\extd E_a.\]
\end{corol} \proof We have to show vanishing of (\ref{cotorG}) for
the coframing $E^{*a}=\eta^{ba}E_a$. However, because $\eta$ is
$\Ad$-invariant and constant, this reduces in terms of $E$ to
vanishing of \[ D_A E_a=\extd E_a+\sum_b (E_{bab^{-1}}-E_a)\wedge
A_b.\] Note that this is a different equation from the torsion
equation solved above. However, since every element of $\CC$ has
order 2, the inverse is irrelevant and the equation then differs
only by a reversal of the $\wedge$. Looking at the equations
solved for zero torsion above, we see that they are invariant
under such a reversal provided we write $A_a=E_b A'_a{}^b$ with
coefficients $A'{}_a{}^b$ from the right. Next we consider the
intersection of the moduli of torion-free and cotorsion-free
connections. Given the bimodule structure, if
$A_u=E_u(\alpha'+1)+E_v\gamma'+E_{uvu}\beta'$, etc., is also
torsion free, we need $R_u(\alpha')=\alpha$, $R_v(\gamma')=\gamma$
and $R_{uvu}(\beta')=\beta$, and similarly for $A_v,A_{uvu}$. As a
result, $R_a(\alpha')=\alpha$ for all $a$, hence $\alpha$ is a
multiple of the identity function (a number) and $\alpha=\alpha'$.
Similarly for $\beta,\gamma$. Finally in this moduli of
torsion-free and cotorsion-free connections we look for regular
connections, i.e. those for which \eqn{regS3}{ A_u\wedge
A_v+A_{uvu}\wedge A_u+A_v\wedge A_{uvu}=0,\quad A_v\wedge
A_u+A_{uvu}\wedge A_v+A_u\wedge A_{uvu}=0,} corresponding to
products of elements from $\CC$ with values $uv$ or $vu$. As
before, we take the first equation, write
$A_u=(\alpha+1)E_u+\gamma E_v+\beta E_{uvu}$ etc., (as found
above), and write all products in terms of our chosen four
2-forms. The coefficients of $E_u\wedge E_v$ and $E_v\wedge E_u$
each yield $\alpha=\gamma$, while those of $E_u\wedge E_{uvu}$ and
$E_{uvu}\wedge E_u$ each yield $\alpha=\beta$. The second equation
above follows in an identical manner and can only give the same
constraints by a symmetry in which we reverse the $\wedge$. Hence
there is a unique regular connection among torsion free and
cotorsion free ones. We write it in the way shown in terms of the
Maurer-Cartan form and $\theta$.

Finally, for any regular connection in our example, the curvature
has to take the form \eqn{regcurvS3}{ F_a=\extd A_a-\sum_b
(A_b\wedge A_a+A_a\wedge A_b)} because the product of all distinct
elements of the conjugacy class lie outside it, so there is no
$A_c\wedge A_d$ term in (\ref{curvG}). For our connections the
second term vanishes since $\sum_b A_b=0$. Also, $\extd\theta=0$
when we put in the values of each $\extd E_a$ and average and use
(\ref{rel2G}). Hence $F_a=\extd E_a$, which is certainly non-zero,
being equal to the quadratic parts in the Maurer-Cartan equation.
\eproof

The explicit $\nabla$ from the general formulae in Section~5.2 is
\[ \nabla E_u=-E_u\tens_H E_u-E_v\tens_H E_{uvu}-E_{uvu}\tens_H
E_v+{1\over 3}\theta\tens_H\theta\]
\[\nabla E_v=-E_v\tens_H E_v-E_u\tens_H
E_{uvu}-E_{uvu}\tens_H E_u+{1\over 3}\theta\tens_H\theta\]
\eqn{nablaS3}{ \nabla E_{uvu}=-E_{uvu}\tens_H E_{uvu}-E_v\tens_H
E_u-E_u\tens_H E_v +{1\over 3} \theta\tens_H\theta} and one may
then verify that indeed torsion and cotorsion vanish as
\[ \nabla\wedge E_a=\extd E_a,\quad
(\nabla\wedge\id-\id\wedge\nabla)(\sum_a E_a\tens_H E_a)=0.\] On
the other hand the similar computation to the latter gives
\[ \nabla(\sum_a E_a\tens_H E_a)=2\sum_{{\rm not}\ a=b=c}
E_a\tens_H E_b\tens_H E_c-2\sum_{\sigma\in
S_3}E_{\sigma(u)}\tens_H E_{\sigma(v)}\tens_H E_{\sigma(uvu)}\ne
0,\] where we keep the left output of $\nabla$ to the far left and
act as a derivation. This is manifestly nonzero (as well as
somewhat basis dependent i.e. not really a computation on
$E_a\tens_H E_a$). Therefore full metric compatibility in the
naive sense does not hold even for this simplest nontrivial
example. This justifies our weaker notion of vanishing cotorsion
as the appropriate generalisation for noncommutative geometry.

We are then able from the general theory above to compute the
Riemann and Ricci curvatures etc., for the Levi-Civita connection
on $S_3$, the latter with respect to a choice of `lift'. One
choice (\ref{worliftG}) is clearly
\[ i(E_u\wedge E_v)=E_u\tens_H E_v-E_w\tens_H E_u,\quad i(E_{uvu}\wedge
E_u)=E_u\tens_H E_{uvu}-E_v\tens_H E_u\] \eqn{worliftS3}{
i(E_v\wedge E_u)=E_v\tens_H E_u-E_{uvu}\tens_H E_v,\quad
i(E_u\wedge E_{uvu})=E_u\tens_H E_{uvu}-E_v\tens_H E_u.} For the
second choice, the basis of $P_{uv}$ and $P_{vu}$ are easily seen
to be the unique vector $\lambda_u=\lambda_v=\lambda_{uvu}=1$ so
that the lift in Proposition~5.3 is \eqn{liftS3}{ i(E_a\wedge
E_b)=E_a\tens_H E_b-{1\over 3}\sum_{cd=ab}E_c\tens_H E_d,\quad
\forall a\ne b.}

\begin{propos} The unique Levi-Civita connection on $S_3$ constructed above
has constant curvature with respect to either of the above two
lifts, with
 \[ {\rm Ricci}=\mu (-g+\theta\tens_H\theta),\]
where $g$ is the metric (\ref{getaS3}) induced by the Killing form
and $\mu=1,2/3$ respectively.
\end{propos}
\proof This is a direct computation from (\ref{ricciG}). The
Riemann tensor is
\[ R E_u= \extd E_u\tens_H E_u+\extd E_v\tens_H E_{uvu}+
\extd E_{uvu}\tens_H E_v,\quad  R E_v=\extd E_v\tens_H E_v+\extd
E_u\tens_H E_{uvu}+\extd E_{uvu}\tens_H E_u \] \eqn{riemS3}{R
E_{uvu}=\extd E_{uvu}\tens_H E_{uvu}+\extd E_v\tens_H E_u +\extd
E_u\tens_H E_v} since $\sum_a \extd E_a=0$. We lift each term by
applying the chosen $i$, then pick out the coefficient of
$E_u\tens$ in $R E_u$ etc., for the trace. \eproof

Thus $S_3$ with its natural Riemannian structure is more or less
an `Einstein space'. We could take
$g_\lambda=g-\lambda\theta\tens_H\theta$ as the metric from the
start without changing anything above (although $\lambda=1$ itself
is degenerate). The scalar curvature itself is the further
contraction of this with the inverse metric. One can similarly
consider several other lifts with the same conclusions but a
different value of $\mu$. Note also that our trace conventions for
Ricci in the classical case would become the first and third
indices of the Riemann tensor, so that we have an opposite sign
convention to the usual one. Hence $S_3$ above for the natural
choices of lift looks more like a compact manifold with constant
positive curvature in usual terms.

Finally, to fix a Dirac operator for the sake of discussion we
choose  the tautological $\gamma$ defined as in (\ref{cangammaG})
by the two-dimensional representation \eqn{rhoWS3}{
\rho_W(u)=\left(\begin{matrix}0&1\\ 1&0\end{matrix}\right),\quad
\rho_W(v)=\left(\begin{matrix}1&0\\-1&-1\end{matrix}\right).} The
braided-Casimir is \eqn{CasS3}{ C={1\over
3}((u-e)^2+(v-e)^2+(uvu-e)^2)=2-{2\over
3}(u+v+uvu),\quad\rho_W(C)=2} so that from (\ref{cangammanorm}),
\eqn{gamnormS3}{ \gamma_a\gamma_b\eta^{ab}=\rho_W(C)=2.} Hence by
Theorem~5.4 (because the elements of $\CC$ all have order 2) the
calculus implied by $\Dsl$ for these will be different from the
one that we have already imposed from quantum group
considerations. In fact $\Dsl$ imposes fewer relations. This is
not a problem from the quantum groups point of view, we can still
use $\Dsl$ perfectly well. The gamma-matrices are explicitly
 \eqn{gammaS3}{ \gamma_u={1\over 3}\left(\begin{matrix}-1&1\\
1&-1\end{matrix}\right), \quad \gamma_v={1\over
3}\left(\begin{matrix}0&0\\ -1&-2\end{matrix}\right), \quad
\gamma_{uvu}={1\over 3}\left(\begin{matrix}-2&-1\\
0&0\end{matrix}\right).} They have some nice identities, however.
In fact for any $\rho_W$ with $\rho_W(uv-e)$ invertible, which is
the case here, one can show by enumeration of the cases and the
identity $\rho_W(e+uv+(uv)^2)=0$ which then holds, that
\eqn{gamrelS3}{\gamma_a\gamma_b+\gamma_b\gamma_a+{2\over
3}(\gamma_a+\gamma_b)={1\over 3}(\delta_{ab}-1),\quad
\sum_a\gamma_a=-1.}

\begin{propos} The Dirac operator on $S_3$ for the above gamma-matrices
and the canonical `Levi-Civita' connection on $S_3$ constructed
above, we have
\[ \Dsl=\del^a\gamma_a-1={1\over 3}\left(\begin{matrix}-\del^u-2\del^{uvu}
-3&\del^u-\del^{uvu}\\
\del^u-\del^v& -\del^u-2\del^v-3\end{matrix}\right).\]
\end{propos}
\proof To find this note first that
$\tau_W^a=\rho_W(a^{-1}-e)=3\gamma_a$ since all elements of $\CC$
have order 2. The canonical connection in terms of components is
$A{}_a{}^b=\delta_a{}^b-{1\over 3}$, hence \[
\Dsl=\del^a\gamma_a-3\sum_a\gamma_a^2+\sum_{a,b}\gamma_a\gamma_b.\]
We then use the gamma-matrix identities above. \eproof

The $-1$ appearing here reflects again a  `constant curvature' now
detected for $S_3$ with its canonical Riemannian structure by the
Dirac operator. Finally we note that while we have focussed here
on the canonical metric induced by the braided-Killing form, one
can similarly consider more general triples $(A,e,e^*)$  and solve
for zero torsion, and zero cotorsion, compute the curvature, etc.
One may then minimise an action defined for example by suitable
contraction of the Ricci curvature, i.e. proceed to finite quantum
gravity. Also, there is no problem introducing Maxwell or
Yang-Mills fields and matter fields since we already have a bundle
formalism, sections etc. This intended application is beyond our
present scope and will be attempted in detail elsewhere. A further
application may be to insert our canonical Dirac operator on $S_3$
into the framework for elementary particle Lagrangians of Connes
and Lott.

\baselineskip 18pt

\end{document}